\documentclass[hidelinks,12pt]{amsart}
\usepackage{bbm}
\usepackage{color}
\usepackage{libertine}
\usepackage{fullpage}
\usepackage{mdframed}
\usepackage{graphicx}
\usepackage{calligra}
\usepackage{wasysym}
\usepackage{csquotes}
\usepackage{bm}
\usepackage{amsmath}
\usepackage[all]{xy}
\usepackage{amssymb}
\usepackage{amsthm}
\usepackage{mathrsfs}
\usepackage[OT2,T1]{fontenc}
\usepackage{centernot}
\usepackage{hyperref}

\newcommand{\DMO}[2]{\DeclareMathOperator{#1}{#2}}

\DMO{\depth}{\mathrm{depth}}

\newtheorem{thm}{Theorem}[section]
\newtheorem{lem}{Lemma}[section]
\newtheorem{coro}{Corollary}[section]
\newtheorem{prop}{Proposition}[section]

\theoremstyle{definition}
\newtheorem{defn}{Definition}[section]
\newtheorem{rmk}{Remark}[section]
\newtheorem{exam}{Example}[section]

\newcommand{\xrar}[1]{\xrightarrow{#1}}

\newcommand{\riso}{\xrar{\sim}}

 \newenvironment{psmat}
  {\left(\begin{smallmatrix}}
  {\end{smallmatrix}\right)}

 \newenvironment{psmatrix}
  {\left(\begin{smallmatrix}}
  {\end{smallmatrix}\right)}

\newcommand{\textboxnobrace}[2] {\parbox{#1em}{\begin{center}{#2}\end{center}}}

\newcommand{\wt}{\widetilde}
\newcommand{\wh}{\widehat}

\newcommand{\ov}{\overline}

\newcommand{\rar}{\rightarrow}

\newcommand{\Rar}{\Rightarrow}

\newcommand{\ncom}[1]{\newcommand{#1}}
\ncom{\sbuset}{\subset}

\newcommand{\hrar}{\hookrightarrow}
\newcommand{\thrar}{\twoheadrightarrow}

\newcommand{\nm}{\mathrm{nm}}

\newcommand{\bs}{\backslash}

\DMO{\bmr}{\mathbbm{r}}
\DMO{\bmf}{\mathbbm{f}}
\DMO{\bmx}{\mathbbm{x}}

\usepackage{enumitem}

\newcommand{\bfone}{\mathbf{1}}
\newcommand{\bA}{\mathbb{A}}

\newcommand{\bC}{\mathbb{C}}

\newcommand{\bG}{\mathbb{G}}

\newcommand{\bN}{\mathbb{N}}

\newcommand{\bP}{\mathbb{P}}
\newcommand{\bQ}{\mathbb{Q}}
\newcommand{\bR}{\mathbb{R}}

\newcommand{\bZ}{\mathbb{Z}}

\newcommand{\bfG}{\mathbf{G}}
\newcommand{\bfH}{\mathbf{H}}

\newcommand{\bfM}{\mathbf{M}}
\newcommand{\bfN}{\mathbf{N}}

\newcommand{\bfP}{\mathbf{P}}

\newcommand{\bfS}{\mathbf{S}}
\newcommand{\bfT}{\mathbf{T}}

\newcommand{\bfY}{\mathbf{Y}}
\newcommand{\bfZ}{\mathbf{Z}}

\newcommand{\cH}{\mathcal{H}}

\newcommand{\cO}{\mathcal{O}}

\DMO{\GWh}{GWh}
\DMO{\DWh}{DWh}

\newcommand{\fg}{\mathfrak{g}}
\newcommand{\fh}{\mathfrak{h}}

\newcommand{\fn}{\mathfrak{n}}

\newcommand{\fz}{\mathfrak{z}}

\newcommand{\fsl}{\mathfrak{sl}}

\DMO{\str}{str}
\DMO{\rel}{-rel}
\DMO{\irrel}{-irrel}
\DMO{\inh}{-inh}
\DMO{\polar}{polar}
\DMO{\Crys}{Crys}
\DMO{\Map}{Map}
\DMO{\Dol}{Dol}
\DMO{\Vect}{Vect}
\DMO{\PGSp}{PGSp}
\DMO{\gon}{gon}
\DMO{\Out}{Out}
\DMO{\comp}{comp}
\DMO{\Inv}{Inv}
\DMO{\Glob}{Glob}
\DMO{\MIC}{MIC}
\DMO{\Isoc}{Isoc}
\DMO{\la}{la}
\DMO{\corank}{corank}
\DMO{\Td}{Td}
\DMO{\hol}{hol}
\DMO{\mir}{mir}
\DMO{\gen}{gen}
\DMO{\BK}{BK}
\DMO{\FL}{FL}
\DMO{\Ann}{Ann}
\DMO{\std}{std}
\DMO{\antidiag}{antidiag}
\DMO{\locadm}{loc.adm}
\DMO{\Inj}{Inj}
\DMO{\LL}{LL}
\DMO{\Dmod}{\emph{D }-mod}
\DMO{\univ}{univ}
\DMO{\Fitt}{Fitt}
\DMO{\WD}{WD}
\DMO{\geom}{geom}
\DMO{\Fl}{Fl}
\DMO{\grad}{grad}
\DMO{\labmda}{\lambda}
\DMO{\Iw}{Iw}
\DMO{\tor}{tor}
\DMO{\coh}{coh}
\DMO{\vol}{vol}
\DMO{\semsim}{ss}
\DMO{\free}{free}
\DMO{\Alg}{Alg}
\DMO{\oth}{otherwise}
\DMO{\Ber}{Ber}
\DMO{\Diff}{Diff}
\DMO{\br}{br}
\DMO{\Isot}{Isot}
\DMO{\prim}{prim}
\DMO{\RAH}{RAH}
\DMO{\Sets}{Sets}
\DMO{\cone}{cone}
\DMO{\Grps}{Grps}
\DMO{\Dec}{Dec}
\DMO{\Flat}{Flat}
\DMO{\AbGps}{AbGps}
\DMO{\Sch}{Sch}
\DMO{\AH}{AH}
\DMO{\cl}{cl}
\DMO{\sk}{sk}
\DMO{\HC}{HC}
\DMO{\cosk}{sk}
\DMO{\ur}{ur}
\DMO{\LocSys}{LocSys}
\DMO{\rk}{rk}
\DMO{\NT}{NT}
\DMO{\cork}{cork}
\DMO{\KS}{KS}
\DMO{\tw}{tw}
\DMO{\MU}{MU}
\DMO{\der}{der}
\DMO{\Art}{Art}
\DMO{\Proj}{Proj}
\DMO{\End}{End}
\DMO{\Betti}{Betti}
\DMO{\Sym}{Sym}
\DMO{\cInd}{cInd}
\DMO{\GL}{GL}
\DMO{\Gal}{Gal}
\DMO{\Br}{Br}
\DMO{\Der}{Der}
\DMO{\Sp}{Sp}
\DMO{\Tan}{Tan}
\DMO{\Spin}{Spin}
\DMO{\Var}{Var}
\DMO{\Nrd}{Nrd}
\DMO{\cusp}{cusp}
\DMO{\Mat}{Mat}
\DMO{\Isom}{Isom}
\DMO{\Stab}{Stab}
\DMO{\SO}{SO}
\DMO{\Res}{Res}
\DMO{\Lie}{Lie}
\DMO{\SU}{SU}
\DMO{\Ad}{Ad}
\DMO{\ad}{ad}
\DMO{\im}{im}
\DMO{\Frob}{Frob}
\DMO{\Fr}{Fr}
\DMO{\red}{red}
\DMO{\an}{an}
\DMO{\Pic}{Pic}
\DMO{\Tor}{Tor}
\DMO{\Hdg}{Hdg}
\DMO{\id}{id}
\DMO{\pr}{pr}
\DMO{\Mor}{Mor}
\DMO{\Ext}{Ext}
\DMO{\ML}{ML}
\DMO{\PGL}{PGL}
\DMO{\SL}{SL}
\DMO{\GU}{GU}
\DMO{\GSp}{GSp}
\DMO{\GSL}{GSL}
\DMO{\Aff}{Aff}
\DMO{\NS}{NS}
\DMO{\gr}{gr}
\DMO{\Ch}{Ch}
\DMO{\QCoh}{QCoh}
\DMO{\Coh}{Coh}
\DMO{\inv}{inv}
\DMO{\Gr}{Gr}
\DMO{\Bun}{Bun}
\DMO{\Hk}{Hk}
\DMO{\GH}{GH}
\DMO{\HT}{HT}
\DMO{\LT}{LT}
\DMO{\Int}{Int}
\DMO{\UU}{U}
\DMO{\OO}{O}
\DMO{\Loc}{Loc}
\DMO{\Conn}{Conn}
\DMO{\sing}{sing}
\DMO{\si}{si}
\DMO{\Sen}{Sen}
\DMO{\MaxSpec}{MaxSpec}
\DMO{\ran}{ran}
\DMO{\coker}{coker}
\DMO{\DIV}{div}
\DMO{\Cl}{Cl}
\DMO{\Frac}{Frac}
\DMO{\VEC}{Vec}
\DMO{\Weil}{Weil}
\DMO{\SPLIT}{split}
\DMO{\Tr}{Tr}
\DMO{\val}{val}
\DMO{\pv}{p.v.}
\DMO{\disc}{disc}
\DMO{\trdeg}{tr.deg}
\DMO{\rad}{rad}
\DMO{\codim}{codim}
\DMO{\dist}{dist}
\DMO{\length}{length}
\DMO{\diam}{diam}
\DMO{\Supp}{Supp}
\DMO{\Ass}{Ass}
\DMO{\ord}{ord}
\DMO{\RE}{Re}
\DMO{\Sh}{Sh}
\DMO{\IM}{Im}
\DMO{\Tot}{Tot}
\DMO{\Bl}{Bl}
\DMO{\lcm}{lcm}
\DMO{\ann}{ann}
\DMO{\arcsinh}{arcsinh}
\DMO{\CHAR}{char}
\DMO{\MOD}{mod}
\DMO{\BB}{BB}
\DMO{\new}{new}
\DMO{\alg}{alg}
\DMO{\Irr}{Irr}
\DMO{\res}{res}
\DMO{\rank}{rank}
\DMO{\naive}{naive}
\DMO{\tors}{tors}
\DMO{\Perf}{Perf}
\DMO{\Sht}{Sht}
\DMO{\Perv}{Perv}
\DMO{\soc}{soc}
\DMO{\Mod}{Mod}
\DMO{\cyc}{cyc}
\DMO{\SC}{sc}
\DMO{\SP}{sp}
\DMO{\Deck}{Deck}
\DMO{\PSL}{PSL}
\DMO{\Area}{Area}
\DMO{\Cont}{Cont}
\DMO{\sgn}{sgn}
\DMO{\Cat}{Cat}
\DMO{\Cov}{Cov}
\DMO{\rig}{rig}
\DMO{\FSch}{FSch}
\DMO{\Rig}{Rig}
\DMO{\Spv}{Spv}
\DMO{\Spa}{Spa}
\DMO{\trace}{trace}
\DMO{\cont}{cont}
\DMO{\aff}{aff}
\DMO{\cor}{cor}
\DMO{\CH}{CH}
\DMO{\Spec}{Spec}
\DMO{\rec}{rec}
\DMO{\LGC}{LGC}
\DMO{\un}{un}
\DMO{\Eval}{Eval}
\DMO{\JH}{JH}
\DMO{\can}{can}
\DMO{\Fss}{Fss}
\DMO{\Speh}{Speh}
\DMO{\Ind}{Ind}
\DMO{\ch}{ch}
\DMO{\nr}{nr}
\DMO{\Swan}{Swan}
\DMO{\St}{St}
\DMO{\Ho}{Ho}
\DMO{\HH}{HH}
\DMO{\trop}{trop}
\DMO{\Jac}{Jac}
\DMO{\vir}{vir}
\DMO{\coll}{coll}
\DMO{\reg}{reg}
\DMO{\dlog}{dlog}
\DMO{\Div}{Div}
\DMO{\ab}{ab}
\DMO{\Tam}{Tam}
\DMO{\Ran}{Ran}
\DMO{\IC}{IC}
\DMO{\Sat}{Sat}
\DMO{\Rat}{Rat}
\DMO{\loc}{loc}
\DMO{\ev}{ev}
\DMO{\st}{st}
\DMO{\pst}{pst}
\DMO{\Fil}{Fil}
\DMO{\cris}{cris}
\DMO{\dR}{dR}
\DMO{\Rep}{Rep}
\DMO{\Sel}{Sel}
\DMO{\spec}{spec}
\DMO{\Spf}{Spf}
\DMO{\JL}{JL}
\DMO{\BGL}{BGL}
\DMO{\Arc}{Arc}
\DMO{\MHS}{MHS}
\DMO{\Nm}{Nm}
\DMO{\holim}{holim}
\DMO{\nInd}{nInd}
\DMO{\sSets}{s\textbf{Sets}}
\DMO{\sArt}{s\textbf{Art}}
\DMO{\BDJ}{BDJ}
\DMO{\GV}{GV}
\DMO{\BM}{BM}
\DMO{\Ord}{Ord}
\DMO{\mult}{mult}
\DMO{\WDRep}{WDRep}
\DMO{\Aut}{Aut}
\DMO{\Hom}{Hom}
\DMO{\sph}{sph}
\DMO{\Def}{Def}
\DMO{\GO}{GO}
\DMO{\diag}{diag}
\DMO{\cond}{cond}
\DMO{\ind}{ind}
\DMO{\irr}{irr}
\DMO{\RHom}{RHom}
\DMO{\sm}{sm}
\DMO{\sss}{ss}
\DMO{\sHom}{sHom}
\DMO{\Tran}{Tran}
\DMO{\Rees}{Rees}
\DMO{\lcv}{lcv} 
\DMO{\SN}{SN}
\DMO{\triv}{triv}
\DMO{\height}{ht}
\DMO{\proj}{proj}
\DMO{\Fun}{Fun}
\DMO{\cts}{cts}
\DMO{\Obj}{Obj}
\DMO{\Sing}{Sing}
\DMO{\Pro}{Pro}
\DMO{\Ig}{Ig}
\DMO{\Ha}{Ha}
\DMO{\BC}{BC}
\DMO{\RZ}{RZ}
\DMO{\supp}{supp}
\DMO{\projdim}{proj.dim}
\DMO{\Zar}{Zar}
\DMO{\Ban}{Ban}
\DMO{\LA}{LA}
\DMO{\ess}{ess}
\DMO{\op}{op}
\DMO{\Func}{Func}
\DMO{\Born}{Born}
\DMO{\Comm}{Comm}
\DMO{\Dr}{Dr}
\DMO{\LC}{LC}
\DMO{\nind}{n-ind}
\DMO{\perf}{perf}
\DMO{\Ob}{Ob}

\usepackage{graphicx} 
\newcommand\smalloplus[1][0.5]{%
   \mathrel{\vcenter{\hbox{\hspace*{0.25em}\scalebox{#1}{$\oplus$}}}}}
\newcommand\smalloplust[1][0.5]{%
   \mathrel{\vcenter{\hbox{\hspace*{0.15em}\scalebox{#1}{$\oplus$}}}}}
\makeatletter
\DeclareRobustCommand{\directsummand}{%
  \mathbin{\mathpalette\o@plus@subset\relax}%
}
\newcommand{\o@plus@subset}[2]{%
  \ooalign{$\m@th#1\smalloplus$\cr$\m@th#1\subset$\cr}%
}
\makeatother

\makeatletter
\DeclareRobustCommand{\directsummandtwo}{%
  \mathbin{\mathpalette\o@plus@subsett\relax}%
}
\newcommand{\o@plus@subsett}[2]{%
  \ooalign{$\m@th#1\smalloplust$\cr$\m@th#1\supset$\cr}%
}
\makeatother

\begin{document}
\title{Generalized Whittaker models beyond $\fsl_{2}$-triples}
\author{Gyujin Oh}\address{Department of Mathematics, Columbia University, 2990 Broadway, New York, NY 10027}
\maketitle
\begin{abstract}
We extend the notion of {generalized Whittaker models} by allowing them to be built upon smooth irreducible  representations of unipotent subgroups of a $p$-adic reductive group that are \emph{not necessarily characters, nor induced from Weil representations}. This notion generalizes the usual notion of generalized Whittaker models, which include the Bessel and  Fourier--Jacobi models. Using Kirillov's orbit method for unipotent groups, we define and prove basic properties of the corresponding Jacquet functors. 
We also provide new instances  of generalized Whittaker models of multiplicity one.
\end{abstract}
\tableofcontents

\section{{Introduction}}

A central problem in the theory of automorphic forms is to determine when a period of an automorphic form vanishes. The local analogue of this problem, the \emph{branching problem}, is as follows: given algebraic groups $\bfH\subset\bfG$ over a local field $F$, describe the restriction of a smooth representation $\pi$ of $\bfG(F)$ to its subgroup $\bfH(F)$. If an irreducible smooth representation of $\bfH(F)$ appears only once as a quotient of $\pi\rvert_{\bfH(F)}$, the Frobenius reciprocity gives a \emph{model} of $\pi$, which means that one could realize $\pi$ as a subrepresentation of a rather concrete representation often realized as the space of functions. 

There are a handful of named \emph{models}, such as the Whittaker models, the Fourier--Jacobi models, and the Bessel models. These models are all of the following form. Let $\bfG$ be a reductive group, and $\bfH=\bfS\bfN$, where $\bfN$ is  the unipotent radical of a parabolic subgroup $\bfP$, and $\bfS\le\bfM$ is a reductive subgroup of the Levi part of $\bfP$. We take an irreducible smooth representation of $\bfH(F)$ of the form $\sigma\otimes\rho$, where $\sigma$ is an irreducible smooth representation of $\bfS(F)$, and $\rho$ is either a character or a Weil representation of $\bfN(F)$, extended over to $\bfH(F)$. In many cases, $\Ind_{\bfH(F)}^{\bfG(F)}(\sigma\otimes\rho)$ is shown to have multiplicity one, including the cases of the aforementioned models, which had many important consequences  in representation theory and the theory of automorphic forms.

Unlike $\sigma$, which can be any irreducible smooth representation of the reductive part $\bfS(F)$, the kinds of representations used in the unipotent part $\bfN(F)$ have been very limited, and in particular ``small'' (in terms of the Gelfand--Kirillov dimension, for example). There are a few reasons for this specification. Firstly, a Weil representation (more so a character) is very rigid, in that it is characterized by its central character. This property yields various ways of explicitly realizing a Weil representation, and is useful in proving various properties of the said models. Another reason is that the models built from a nilpotent orbit (or an $\fsl_{2}$-triple) in $\Lie\bfG$ have a relationship with the local character expansion of a representation of $\bfG(F)$ (e.g., \cite{MoeglinWaldspurger}), and either a character or a Weil representation can only appear at the unipotent part for all such models.

In this paper, we study a more general kind of model, where we allow $\rho$ to be \emph{any smooth irreducible representation} of $\bfN(F)$. We refer to such models as the \emph{generalized Whittaker models}\footnote{In the literature, this name  referred only to   the models  coming from $\fsl_{2}$-triples (hence the ``beyond $\fsl_{2}$-triples'' in the title); see Remark~\ref{rem:history} for a comparison between our definition and those in the literature.}. A key tool in studying these models is \emph{Kirillov's orbit method} for unipotent groups \cite{Kirillov,Moore} (reviewed in \S\ref{sec:Kirillov}), which gives a complete  understanding of smooth  irreducible representations of $\bfN(F)$ in terms of the coadjoint orbits of $\bfN(F)$. For example, one can explicitly realize any smooth irreducible representation of $\bfN(F)$ on a space of functions, and one has a combinatorial recipe to tell when two smooth irreducible representations of $\bfN(F)$ are isomorphic to each other (albeit much more complicated than the Stone--von Neumann theorem in the case of Weil representations). In this paper, we prove several basic properties of such models, such as the existence of a metaplectic cover of finite degree in great generality (see \S\ref{sec:metaplectic}). We also define the appropriate generalizations of the \emph{Jacquet functors} and the \emph{degenerate Whittaker models}, and prove fundamental properties such as their exactness (see \S\ref{sec:GWh}). Using these, we show that there are genuinely new  instances of generalized Whittaker models of multiplicity one:
\begin{thm}[Rough form]
There exist \emph{high-depth} generalized Whittaker models of multiplicity one for the minimal parabolic subgroup of $\GL_{n}$ (for $n\ge4$) or $\Sp_{2n}$ (for $n\ge3$).
\end{thm}
See \S\ref{sec:exam} for the precise statements. Here, we call a generalized Whittaker model \emph{high-depth} if the representation $\rho$ of the unipotent part is neither a character nor induced from a Weil representation. By definition, any prior instances of generalized Whittaker models are not high-depth (referred to as \emph{low-depth} in this paper).

This new notion of generalized Whittaker models naturally arises in the problem of \emph{Fourier expansions} of automorphic forms. For example, suppose that $\bfG$ is a reductive algebraic group over $\bQ$, $\bfP$ is a parabolic subgroup, and $\bfN$ is its unipotent radical. Over its center $\bfZ\le\bfN$, we can take a Fourier expansion of a cuspidal automorphic form $f$ on $\bfG(\bA)$:
\[f(g)=\sum_{\xi\in\wh{\bfZ(\bA)}} f_{\xi}(g),\quad f_{\xi}(g):=\int_{\bfZ(\bQ)\bs \bfZ(\bA)}f(zg)\xi(z)dz.\]
If we choose a cuspidal automorphic representation $\pi$ of $\bfG(\bA)$ and vary $f$ within $\pi$, then, for a chosen character $\xi\in\wh{\bfZ(\bA)}$, $f\mapsto f_{\xi}$ has a property that $(z\cdot f)_{\xi}=\xi^{-1}(z)f_{\xi}$ for any $z\in \bfZ(\bQ_{p})$ for any prime $p$. If $\bfN$ is either abelian or a Heisenberg group, then this would imply that $f\mapsto f_{\xi}$ transforms as a character or a Weil representation of $\bfN(\bQ_{p})$ by the Stone--von Neumann theorem. However, in general, the central character is not enough to determine a representation of $\bfN(\bQ_{p})$. As a concrete example, when $\bfG=\GL_{4}$ and $\bfP$ is the standard upper triangular Borel subgroup, the coadjoint orbits of $\bfN(\bQ_p)$ are calculated in Example~\ref{exam:GL4}. The calculation suggests that, for $\xi\ne1$, $f_{\xi}$ should be further decomposed into
\[f_{\xi}(g)=\sum_{\zeta\in\wh{\bfY(\bA)}}f_{\xi,\zeta}(g),\quad f_{\xi,\zeta}(g):=\int_{\bfY(\bQ)\bs\bfY(\bA)}f(yg)\zeta(y)dy,\quad \bfY:=\left\lbrace\begin{pmatrix}1&0&0&0\\0&1&*&0\\0&0&1&0\\0&0&0&1\end{pmatrix}\right\rbrace.\]
The refined Fourier coefficient $f_{\xi,\zeta}$ is supported in a \emph{degenerate Whittaker model} in the sense of Definition~\ref{def:DWh}, and, further decomposing $f_{\xi,\zeta}$ based on the stabilizer group action, we obtain functionals supported in a generalized Whittaker model. Note that there has recently been an exciting development of the algebraicity  of Fourier coefficients of automorphic forms of more general types (e.g., \cite{Pollack1,Pollack2,Pollack3}); we hope that the generalized Whittaker models in this paper's sense will be found useful in definining the notion of algebraicity for more general types of automorphic forms. See also \S\ref{sec:Fourier}.

Another interesting feature of the generalized Whittaker models in our sense comes from the ``flexibility'' of $\bfN(F)$-coadjoint orbits. Namely, there are continuous families of coadjoint orbits of $\bfN(F)$ of the same dimension, as opposed to $\bfG(F)$. In particular, the generalized Whittaker model depends not only on the size of the $\bfN(F)$-coadjoint orbit corresponding to the representation of $\bfN(F)$, but also its $\bfM(F)$-orbit in the space of $\bfN(F)$-coadjoint orbits. Guided by the example of $\GL_{4}$, we expect that  the generalized Whittaker model of an $N$-coadjoint orbit is related to that of another $N$-coadjoint orbit obtained as a deformation of the original $N$-coadjoint orbit. See \S\ref{sec:GL4} for more details.

In \S\ref{sec:BZSV}, we associate a Hamiltonian space associated to  the generalized Whittaker model in the sense of this paper, and explain where it sits in the optic of the recent program of relative Langlands duality of \cite{BZSV}. In particular, we show that the Hamiltonian space does not satisfy certain  axioms such as affineness and neutrality required by \emph{op. cit.} Finally, in Appendix, we collect some lemmas in smooth representation theory.
\subsection*{Acknowledgements} To be added later.

\smallskip

\noindent\textit{Notations and Conventions.} Let $p$ be {an odd  prime}, and $F$ be a finite extension of $\bQ_{p}$, with a uniformizer $\pi$. Let $S^{1}\subset\bC^{\times}$ be the topological multiplcative group of norm $1$ complex numbers. Let $\varepsilon\colon F\rar S^{1}$ be the normalized unitary character of $F$ {in the sense of Tate's thesis \cite{TateThesis}.} 

Let $\bfN$ be a unipotent algebraic group over $F$ (i.e. there exists a composition series whose successive quotients are $\bG_{a}$, see \cite[Expos\'e XVII, D\'efinition 1.3, Corollaire 4.1.3]{SGA3-2}), and let $N=\bfN(F)$. Let $\Lie\bfN$ be the Lie algebra of $\bfN$, which is an affine algebraic variety over $F$. Then, there is an isomorphism of $F$-varieties $\exp:\Lie\bfN\rar\bfN$ (note that the power series defining the exponential gets truncated to a polynomial in this case). This induces an isomorphism $\exp:\Lie N\rar N$, where $\Lie N:=(\Lie\bfN)(F)$. We denote the inverse of $\exp$ as $\log$. We also note that $N$ is unimodular. 

A \emph{td-group} is a locally profinite group (i.e. a Hausdorff totally disconnected topological group where any open neighborhood of $1$ contains an open compact subgroup) whose quotient space modulo an open compact subgroup is countable. A typical example of a td-group is $\bfG(F)$ for a linear algebraic group $\bfG$ over $F$. Every representation in this paper has coefficients over $\bC$. A smooth representation of a td-group is a representation where the stabilizer of a vector is an open is open. For a td-group $G$, let $\Rep^{\sm}(G)$ be the category of smooth representations of $G$. Note that Schur's lemma holds (e.g., \cite[1.2.6]{BushnellHenniart}). Let $\wh{G}$ be the set of all irreducible smooth representations of $G$. When $G$ is a unipotent group over $F$, $\wh{G}$ is parametrized by the coadjoint orbits of $G$ by Kirillov's orbit method, Theorem~\ref{thm:Kirillov}. In this context, whenever something is indexed by an element of $\wh{G}$, we will use $\psi$, a linear functional on $\Lie G$, and $\rho_{\psi}$, its corresponding irreducible smooth representation of $G$, interchangeably.

Given any $G$-representation $\pi$, let $\pi_{\infty}$ be the space of smooth vectors in $\pi$, which forms a smooth representation of $G$. Note that $\Rep^{\sm}(G)$ is a $\bC$-linear abelian category that has enough injectives, so we may talk about the right derived functor of $\Hom_{G}$ (the space of $G$-linear maps), which we denote $\Ext^{i}_{G}$. In $\Rep^{\sm}(G)$, we have the involution called the \emph{smooth dual} (or the \emph{contragredient}) $\pi\mapsto\pi^{\vee}:=(\pi^{*})_{\infty}$ ($\pi^{*}$ is the linear dual). For $\pi,\sigma\in\Ob(\Rep^{\sm}(G))$, $\Hom_{G}(\pi,\sigma^{\vee})=\Hom_{G}(\sigma,\pi^{\vee})$.

For a closed subgroup $H\le G$, let $\delta:G\rar\bR_{>0}$ be the modulus character given by the sheaf of locally constant real-valued measures on $G/H$, and let $\delta^{1/2}$ be the square-root of $\delta$. Using this, we  define the (normalized) induction $\Ind_{H}^{G}:\Rep^{\sm}(H)\rar\Rep^{\sm}(G)$ and the (normalized) compact induction $\ind_{H}^{G}:\Rep^{\sm}(H)\rar\Rep^{\sm}(G)$. By the same proof as \cite[Proposition 29]{Bernstein-Notes}, $\left(\ind_{H}^{G}(\pi)\right)^{\vee}=\Ind_{H}^{G}(\pi^{\vee})$ for $\pi\in\Ob(\Rep^{\sm}(H))$.

We write $Z(G)$ for the center of $G$. The $n\times n$ identity matrix is denoted $I_{n}$. The expression $\diag(a_{1},\cdots,a_{n})$ means the diagonal $n\times n$ matrix with diagonal entries $a_{1},\cdots,a_{n}$ in order from top-left to bottom-right. The group of $n\times n$ matrices with entries in $F$ is denoted $M_{n}(F)$.
\section{{Kirillov's orbit method}}\label{sec:Kirillov}
We first review Kirillov's \emph{orbit method} \cite{Kirillov, Moore}, which gives an explicit parametrization of smooth admissible irreducible representations of $p$-adic unipotent groups. Throughout this section, $N$ denotes the group of $F$-points of a unipotent algebraic group over $F$.

\begin{thm}[{{Kirillov's orbit method}}]\label{thm:Kirillov}
Let $(\Lie N)^{*}$ be the $F$-linear dual of $\Lie N$. Then, there is an explicit bijection
\[
(\Lie N)^{*}/N\riso\wh{N},\quad \cO_\psi\mapsto\rho_{\psi},
\]
where, for $\psi\in(\Lie N)^{*}$, $\cO_{\psi}\in (\Lie N)^{*}/N$ is the coadjoint orbit of $\psi$ (i.e. the $N$-orbit of $\psi$ under the coadjoint action of $N$ on $(\Lie N)^{*}$). More precisely, given $\psi\in(\Lie N)^{*}$, a construction of  an irreducible smooth  representation $\rho_{\psi}$ of $N$ is given as follows.
\begin{enumerate}
\item Choose a \emph{polarization $\fh\subset \Lie N$ subordinate to $\psi$}, i.e., a maximal Lie  subalgebra such that $\psi\rvert_{[\fh,\fh]}=0$.
\item Let $H:=\exp(\fh)\le N$. Consider the unitary character\footnote{This is indeed a character due to the condition $\psi\rvert_{[\fh,\fh]}=0$.} $\chi_{\psi}:H\rar S^{1}$ defined by 
\[\chi_{\psi}(h):=\varepsilon(\psi(\log h)).\]
\item We define $\rho_{\psi}:=\Ind_{H}^{N}\chi_{\psi}\cong\ind_{H}^{N}\chi_{\psi}$.
\end{enumerate}
Furthermore, $\rho_{\psi}$ is unitary and admissible.
\end{thm}
\begin{proof}
This Theorem is proved by an inductive procedure in \cite{Kirillov}, specifically Theorem~5.1 of \emph{op.~cit.}, in the case of simply connected nilpotent Lie groups. As remarked in \cite[Theorem 3]{Moore}, the same proof works verbatim for unipotent groups over a $p$-adic local field of characteristic $0$ (this was explicitly worked out in \cite[\S3]{Matringe}). The  isomorphism $\Ind_{H}^{N}\chi_{\psi}\cong\ind_{H}^{N}\chi_{\psi}$ is a consequence of \cite[Theorem 4]{vanDijk} (which implies that $\Ind_{H}^{N}\chi_{\psi}$ is admissible). The unitarity  is \cite[Proposition 3.1]{vanDijk}.
\end{proof}

\begin{exam}\label{exam:oscillator}
Let $(W,\langle,\rangle)$ be a symplectic space over $F$. Then, the \emph{Heisenberg group} $H(W)$ is the group defined on the set $F\times W$ by the multiplication law
\[(a,v)\cdot (b,w):=\left(a+b+\frac{1}{2}\langle v,w\rangle,v+w\right).\]
Then, $H(W)$ is (the $F$-points of) a unipotent group over $F$, and the center $Z\subset H(W)$ consists of the elements of the form $(a,0)\in F\times W$. Note that $\Lie Z=[\Lie H(W),\Lie H(W)]$. 

Let $\psi\in(\Lie H(W))^{*}$. If $\psi\rvert_{\Lie Z}=0$, then $\cO_\psi$ is a singleton, and $\rho_{\psi}=\chi_{\psi}$ is a character. On the other hand, if $\psi\rvert_{\Lie Z}\ne0$, then one can find $\psi'\in\cO_\psi$  such that, under the vector space decomposition $\Lie H(W)=\Lie Z\oplus W$ coming from $H(W)=F\times W$, $\psi'\rvert_{W}=0$. Then, for a maximal isotropic subspace $Y\subset W$ (i.e. $\langle Y,Y\rangle=0$ and $\dim Y=\frac{1}{2}\dim W$), $\fh:=\Lie Z\oplus Y\subset\Lie H(W)$ is a polarization subordinate to $\psi'$, and $\rho_{\psi}=\Ind_{H}^{H(W)}\chi_{\psi'}$ is (the Schr\"odinger model of) the \emph{Weil representation}  of $H(W)$ with central character $\chi_{\psi'}$.
%
\end{exam}

\begin{coro}\label{coro:dual of Kirillov}
For $\psi\in(\Lie N)^{*}$, $\left(\rho_{\psi}\right)^{\vee}\cong\rho_{-\psi}$.
\end{coro}
\begin{proof}
\[\left(\rho_{\psi}\right)^{\vee}\cong\left(\ind_{H}^{N}\chi_{\psi}\right)^{\vee}\cong\Ind_{H}^{N}\left(\chi_{\psi}^{\vee}\right)\cong\Ind_{H}^{N}\chi_{-\psi}=\rho_{-\psi}.\]
\end{proof}

\begin{defn}
Let $\psi\in(\Lie N)^{*}$. The \emph{depth} of $\rho_{\psi}$ (or depth of $\psi$), denoted $\depth\rho_{\psi}$ (or $\depth\psi$), is defined as
\[\depth\rho_{\psi}:=\min\lbrace n\ge0~:~\rho_{\psi}\rvert_{N_{n+1}}\text{ is trivial}\rbrace,\]
where $N=N_{1}\trianglerighteq N_{2}\trianglerighteq\cdots$ is the lower central series, i.e. $G_{n}:=[G_{n-1},G]$. We say that $\rho_{\psi}$ is \emph{low-depth} (\emph{high-depth}, respectively) if $\depth\rho_{\psi}\le2$ ($\depth\rho_{\psi}\ge3$, respectively). 

We also define the \emph{depth of $N$} to be the smallest $n$ such that $N_{n+1}=0$. For example, if $N$ is abelian (meta-abelian, respectively), $\depth N=1$ ($\depth N=2$, respectively).
\end{defn}

Note that the low-depth representations of $N$ are those which had been intensively studied.

\begin{lem}\label{lem:depth}
Let $\psi\in(\Lie N)^{*}$. 
\begin{enumerate}
\item We have $\depth\psi=1$ if and only if $\rho_{\psi}$ is a nontrivial character of $N$.
\item We have $\depth\psi=2$ if and only if there is a Weil representation $\sigma$  of the Heisenberg group $H(W)$ of a vector space $W$ over $F$ with a symplectic pairing $\langle,\rangle:W\times W\rar F$ (see Example~\ref{exam:oscillator}) and a quotient $f:N\thrar H(W)$  such that $\rho_{\psi}=f^{*}\sigma$.
\end{enumerate}
\end{lem}
\begin{proof}
\begin{enumerate}
\item If $\rho_{\psi}$ is a nontrivial character, then $\rho_{\psi}$ is trivial on $[N,N]$, so $\depth\psi=1$. Conversely, if $\depth\psi=1$, then $\rho_{\psi}$ is an irreducible representation of $N/N_{1}=G^{\ab}$, so $\rho_{\psi}$ must be a character; it is a nontrivial character as it is not trivial on the whole $N$.
\item Suppose that there is a quotient $f:N\rar H(W)$ and a Weil representation $\sigma$ of $H(W)$ such that $\rho_{\psi}=f^{*}\sigma$. As $H(W)_{2}=0$, $\rho_{\psi}$ is trivial on $G_{2}$, so $\depth\psi\le 2$. As $\rho_{\psi}$ is not a character ($\sigma$ is infinite-dimensional), by (1), we know that $\depth\psi=2$. Conversely, suppose that $\depth\psi=2$. By quotienting out by $N_{2}$, we may assume that $N_{2}=0$. As $\depth\psi>1$, $N_{1}\ne0$ and $\rho_{\psi}$ is nontrivial on $N_{1}$. As $N_{1}$ is contained in the center of $N$, it follows that $\rho_{\psi}$ is nontrivial on the center of $N$. By the inductive procedure for the proof of Theorem~\ref{thm:Kirillov}, it follows that $\rho_{\psi}$ is the pullback of a representation $\sigma$ of a quotient $N\thrar N'$ where the center $Z(N')$ of $N'$ is one-dimensional (and $\sigma$ is nontrivial on $Z(N')$).

We claim that $N'$ is the Heisenberg group of a symplectic space over $F$, and $\sigma$ is a Weil representation. Note that $W:=N'/Z(N')$ is abelian and unipotent, so, by \cite[Expos\'e XVII, Corollaire 4.1.3]{SGA3-2}, it is a finite-dimensional vector space over $F$ (and by the same reason, $Z(N')\cong F$). Consider the commutator morphism $[,]:N'\times N'\rar Z(N')$. As $Z(N')$ is in the kernel, it give rise to a morphism $[,]:W\times W\rar Z(N')$. It is easy to see that this is a bilinear morphism of commutative group schemes, and is alternating. Therefore, this gives rise to an alternating pairing $\langle,\rangle:W\times W\rar F$. It is then straightforward to exhibit an isomorphism $N'\cong H(W)$ to the (degenerate) Heisenberg group $H(W)$ associated with the alternating pairing $\langle,\rangle:W\times W\rar F$. Here, the degenerate Heisenberg group means $H(W)=H(W')\times W''$, where $W''=\ker\langle,\rangle$, and $\langle,\rangle\rvert_{W'}$ is non-degenerate (i.e., symplectic). As $Z(N')=Z(H(W))$ is one-dimensional, $H(W)$ must be a Heisenberg group (i.e. $\langle,\rangle$ is non-degenerate). As $\sigma$ is an infinite-dimensional smooth admissible irreducible  representation of $H(W)$, it must be a  Weil representation. 
\end{enumerate}
\end{proof}

\section{Metaplectic covers of finite degree}\label{sec:metaplectic}
In the case of the Weil representation $\rho_{\psi}$ of a Heisenberg group $H(W)$ (Example~\ref{exam:oscillator}), $\Sp(W)$ acts on $H(W)$, which preserves the isomorphism class of $\rho_{\psi}$. This gives rise to a  \emph{projective representation} of $\Sp(W)$ on $\rho_{\psi}$. A priori, the unitarity of $\rho_{\psi}$ only asserts that $\rho_{\psi}$ exists as a linear representation of a central $S^{1}$-cover $\wt{\Sp(W)}$ of $\Sp(W)$, called the \emph{metaplectic $S^{1}$-cover}. However, it is well-known (\cite{Weil}) that $\wt{\Sp(W)}$ reduces to a double cover $\wh{\Sp(W)}$ of $\Sp(W)$ (also called the \emph{metaplectic group}), and $\wh{\Sp(W)}$ acts linearly on $\rho_{\psi}$. 

One may try to study a similar problem with more general irreducible smooth admissible representations constructed by Kirillov's orbit method, Theorem~\ref{thm:Kirillov}. The setup of the rest of the paper is as follows. Let $\bfG$ be a reductive $p$-adic group over $F$, and let $\bfP\le \bfG$ be a parabolic subgroup with the Levi decomposition $\bfP=\bfM\bfN$. Let $G=\bfG(F)$, $P=\bfP(F)$, $M=\bfM(F)$, $N=\bfN(F)$. Let $(\rho,V)$ be a smooth admissible irreducible  representation of $N$. The conjugation action by $M$ gives an action of $M$ on $\wh{N}$. Let $S_{\rho}\le M$ be the identity component of the stabilizer of the isomorphism class $[\rho]\in\wh{N}$ by the said action. By Schur's lemma, this gives rise to a projective representation of $S_{\rho}$ on $V$, i.e. a homomorphism $S_{\rho}\rar\PGL(V)$, and by the unitarity of $\rho$ (see Theorem~\ref{thm:Kirillov}), this further lifts to a homomorphism $S_{\rho}\rar\GL(V)/S^{1}$, where $S^{1}\subset\bC^{\times}$ is the group of norm one complex numbers. 

In this context, is easy to define the analogue of the metaplectic $S^{1}$-cover.
\begin{defn}[Metaplectic $S^{1}$-cover]Let $\wt{S_{\rho}}$ be the topological $S^{1}$-cover of $S_{\rho}$, called the \emph{metaplectic $S^{1}$-cover} of $S_{\rho}$,
\[1\rar S^{1}\rar \wt{S_{\rho}}\rar S_{\rho}\rar1,\]
given by pulling back the short exact sequence
\[1\rar S^{1}\rar \GL(V)\rar \GL(V)/S^{1}\rar1,\]
along the map $S_{\rho}\rar \GL(V)/S^{1}$. By definition, there is a natural homomorphism $\wt{S_{\rho}}\rar\GL(V)$. 
\end{defn} 
However, for the sake of smooth representation theory, it is desirable to work with td-groups, whereas the metaplectic $S^{1}$-cover $\wt{S_{\rho}}$ is not totally disconnected. In this section, we show that the metaplectic $S^{1}$-cover reduces to a finite-degree cover in most cases.
\begin{thm}\label{thm:finite cover exists}
Let $\psi\in(\Lie N)^{*}$, and let $(\rho_{\psi},V)$ be the irreducible smooth admissible  representation of $N$ associated with $\psi$ by Theorem~\ref{thm:Kirillov}. Suppose that Assumption~\emph{(\ref{tag:reductive})} holds. Then, there exists a covering $\wh{S}\rar S_{\rho_{\psi}}$ of finite degree,
\[1\rar\mu_{n}\rar\wh{S}\rar S_{\rho_{\psi}}\rar1,\]
such that the extension $\wt{S}$ is induced from $\wh{S}$ by the inclusion $\mu_{n}\rar S^{1}$.  This finite cover $\wh{S}$ is a td-group, and $V$ is a smooth representation of $\wh{S}$.
\end{thm}
{
The Assumption (\ref{tag:reductive}) is as follows.
\[
\tag{$\blacklozenge$}\label{tag:reductive}\textboxnobrace{36}{The stabilizer $S_{\rho_{\psi}}$ is the set of $F$-points of a linear algebraic group, and if we let $U:=R_{u}(S_{\rho_{\psi}})$ be the unipotent radical and $R:=S_{\rho_{\psi}}/U$ be the reductive quotient,  there exists  $t\in Z(R)$ such that the $\ker(t:\Lie U\rar\Lie U)=0$.}\]
We expect that Assumption~(\ref{tag:reductive}) holds in great generality. For example, the group theoretic condition of Assumption~(\ref{tag:reductive}) holds for any parabolic subgroup of a reductive group (see \cite[\S3, Proposition]{Sury}). In particular, all examples known to the author satisfy Assumption (\ref{tag:reductive}). It would be interesting to either prove Assumption (\ref{tag:reductive}) holds in general, or to find an example where Assumption (\ref{tag:reductive}) {does not hold}. }
\begin{defn}[Metaplectic cover]\label{def:metaplectic}
Let $\psi\in(\Lie N)^{*}$, and let $(\rho_{\psi},V)$ be the irreducible smooth admissible  representation of $N$ associated with $\psi$ by Theorem~\ref{thm:Kirillov}. A \emph{metaplectic cover} is a covering,
\[
1\rar\mu_{n}\rar\wh{S_{\rho_{\psi}}}\rar S_{\rho_{\psi}}\rar1,
\]
of \emph{minimal} degree $n$, such that there exists a smooth (linear) representation of $\wh{S_{\rho_{\psi}}}\ltimes N$ on $V$ such that its restriction to $N$ coincides with $\rho_{\psi}$. We define the \emph{degree} of the metaplectic cover $\wh{S_{\rho_{\psi}}}$ to be $n$. 
\end{defn}
Thus, Theorem~\ref{thm:finite cover exists} should be seen as asserting that Definition~\ref{def:metaplectic} is well-defined in most cases. 

Before we show Theorem~\ref{thm:finite cover exists}, we show some auxiliary results asserting the existence of a metaplectic cover.

\begin{lem}\label{lem:metaplectic splits for characters}
If $\rho$ is a  character of $N$, then $\wh{S_{\rho}}=S_{\rho}$ is the metaplectic cover.
\end{lem}
\begin{proof}
Regarding $\rho$ as a homomorphism $N\rar S^{1}$, that $m\in M$ is contained in $S_{\rho}$ means that $\rho$ is isomorphic to $\rho_{m}:N\rar S^{1}$ defined by $\rho_{m}(n):=\rho(mnm^{-1})$. This however means that $\rho_{m}(n)=\rho(n)$ for all $n\in N$. Therefore, $\rho$ trivially extends to a character of $S_{\rho}\rvert_N$, which implies that $\wt{S_{\rho}}$ splits.
\end{proof}

More generally, the following holds.

\begin{lem}\label{lem:split over parabolic}
Let $\psi\in(\Lie N)^{*}$, and let $\fh\subset\Lie N$ be a polarization subordinate to $\psi$. If $S_{\rho_{\psi}}$ stabilizes the flag $0\subset\fh\subset \Lie N$, then  $\wh{S_{\rho_{\psi}}}=S_{\rho_{\psi}}$ is the metaplectic cover.
\end{lem}
\begin{proof}
We show this by explicitly constructing a genuine action of $S_{\rho_{\psi}}$ on 
\[
V=\Ind_{H}^{N}\chi_{\psi}=\lbrace f:N\rar \bC~\text{uniformly locally constant}~:~f(hn)=\chi_{\psi}(h)f(n)~\forall h\in H,n\in N\rbrace,
\]
as follows (the modulus character does not appear, as $N$ is unimodular). For $g\in S_{\rho_{\psi}}$, let $\varphi_{g}:N\rar N$ be the automorphism induced by conjugation  $n\mapsto gng^{-1}$. By assumption, $\varphi_{g}$ restricts to $\varphi_{g}:H\rar H$, and, importantly, $\chi_{\psi}(\varphi_{g}(h))=\chi_{\psi}(h)$ for $h\in H$. We now define the action of $N_{S_{\rho_{\psi}}}(\fh)$ on $V$ by
\[(g\cdot f)(n):=f(\varphi_{g}(n)).\]
This is well-defined as, for $h\in H$ and $n\in N$,
\[(g\cdot f)(hn)=f(\varphi_{g}(hn))=f(\varphi_{g}(h)\varphi_{g}(n))=\chi_{\psi}(\varphi_{g}(h))f(\varphi_{g}(n))=\chi_{\psi}(h)(g\cdot f)(n).\]
\end{proof}
The following asserts that, for \emph{low-depth} representations of $N$ (i.e. $\depth\le2$), the metaplectic cover exists and is of low degree.
\begin{coro}\label{coro:deg}Let $\psi\in(\Lie N)^{*}$.
\begin{enumerate}
\item If $\depth\psi=1$, then $\deg\wh{S_{\rho_{\psi}}}=1$, i.e. $\wh{S_{\rho_{\psi}}}=S_{\rho_{\psi}}$.
\item If $\depth\psi=2$, then $\deg\wh{S_{\rho_{\psi}}}\le2$.
\end{enumerate}
\end{coro}
\begin{proof}
\begin{enumerate}
\item This is Lemma~\ref{lem:metaplectic splits for characters}.
\item If $\depth\psi=2$, by Lemma~\ref{lem:depth}(2), $\rho_{\psi}=f^{*}\sigma$ for a Weil representation $\sigma$ of the Heisenberg group $H(W)$ of a symplectic vector space $W$ over $F$ along a surjective homomorphism $f:N\thrar H(W)$. Although a general automorphism on $N$ does not descend to an automorphism on $H(W)$ (i.e., $H(W)$ is not a \emph{characteristic quotient} of $N$), we claim that the action of $S_{\rho_{\psi}}$ on $N$ does preserve the quotient $H(W)$, inducing an action of $S_{\rho_{\psi}}$ on $H(W)$. Firstly, $f$ factors through $N/N_{2}$ (recall that $N_{2}=[N,[N,N]]$), which is a characteristic quotient of $N$, so the action of $S_{\rho_{\psi}}$ on $N$ descends to an action on $N/N_{2}$. Note also that  the coadjoint orbit $\cO_{\psi}$ is contained in $(\Lie (N/N_{2}))^{*}\subset(\Lie N)^{*}$. 

The inductive procedure alluded in the proof of Lemma~\ref{lem:depth}(2) shows that $H(W)$ is obtained from $N/N_{2}$ by a chain of quotients 
\[
N/N_{2}=H_{0}\thrar H_{1}\thrar\cdots\thrar H_{n}=H(W),\]
where, $H_{i}$ is meta-abelian for $0\le i\le n$ (i.e., $[H_{i},[H_{i},H_{i}]]=1$), and $\ker(H_{i}\thrar H_{i+1})\subset Z(H_{i})$ for $0\le i<n$. We claim that  $\cO_{\psi}\subset (\Lie H_{n})^{*}$. We already know that $\cO_{\psi}\subset(\Lie H_{0})^{*}$, so it is sufficient to show that $\cO_{\psi}\subset(\Lie H_{i})^{*}$ implies that $\cO_{\psi}\subset(\Lie H_{i+1})^{*}$ for $0\le i<n$. Note that, as $\ker(H_{i}\thrar H_{i+1})$ is central in $H_{i}$, any conjugate of $\psi$ by an element of $H_{i}$ is trivial on $\Lie(\ker(H_{i}\thrar H_{i+1}))$. Therefore, $\cO_{\psi}\subset (\Lie(H_{i+1}))^{*}$, as desired. Moreover, for any $\psi'\in\cO_{\psi}$, $H(W)$ can be recovered from $\psi'$ by following the above procedure. Therefore, it follows that the action of $S_{\rho_{\psi}}$ preserves $\ker(N/N_{2}\thrar H(W))$, so the action of $S_{\rho_{\psi}}$ on $N/N_{2}$ descends to an action of $S_{\rho_{\psi}}$ on $H(W)$. As the natural map $f^{*}:\wh{H(W)}\rar \wh{N}$ is injective, the action of $S_{\rho_{\psi}}$ on $H(W)$ preserves the isomorphism class of the Weil representation $\sigma$.

Note that action of $M$ on $\Lie N$ is $F$-linear. As $\Lie N\thrar \Lie H(W)$ is $F$-linear, it follows that the action of $S_{\rho_{\psi}}$ on $\Lie H(W)$ is $F$-linear. Let $\Aut_{F}(H(W),\sigma)$ be the group of continuous automorphisms of $H(W)$ which are $F$-linear on $\Lie H(W)$ and which preserve the isomorphism class of $\sigma$. Then, the action of $S_{\rho_{\psi}}$ on $H(W)$ gives rise to a group homomorphism $S_{\rho_{\psi}}\rar \Aut_{F}(H(W),\sigma)$. If $\varphi:H(W)\rar H(W)$ is an element of $\Aut_{F}(H(W),\sigma)$, then, using the expression $H(W)=F\times W$ as in Example~\ref{exam:oscillator}, we can write $\varphi(a,v)=(\varphi_{1}(a,v),\varphi_{2}(a,v))$. As the second component $W$ is the abelianization of $H(W)$, it follows that $\varphi_{2}(a,v)$ does not depend on $a$, and is $F$-linear in $v$. Thus, $\varphi_{2}\in\GL(W)$. Using the multiplication law, we have, for any $a,b\in F$ and $v,w\in W$,
\[\varphi_{1}(a,v)+\varphi_{1}(b,w)+\frac{1}{2}\langle\varphi_{2}(v),\varphi_{2}(w)\rangle=\varphi_{1}\left(a+b+\frac{1}{2}\langle v,w\rangle,v+w\right).\]Putting $v=0$, we have
\[\varphi_{1}(a,0)+\varphi_{1}(b,w)=\varphi_{1}(a+b,w).\]
Therefore, for any $a\in F$ and $v\in W$, $\varphi_{1}(a,v)=\varphi_{1}(a,0)+\varphi_{1}(0,v)$. Note that $\varphi_{1}(-,0):F\rar F$ and $\varphi_{1}(0,-):W\rar F$ are both $F$-linear as noted above. Therefore, $\varphi_{1}(a,0)=\lambda a$ for $\lambda\in F$. As the central character of $\sigma$ is preserved under the action of $\varphi$, it follows that $\lambda=1$. Thus, the multiplication law simplifies into
\[\frac{1}{2}\langle\varphi_{2}(v),\varphi_{2}(w)\rangle=\frac{1}{2}\langle v,w\rangle,\]
i.e., $\varphi_{2}\in\Sp(W)$. In particular, we see that
\[\Aut_{F}(H(W),\sigma)=\begin{pmatrix}1& *\\ 0&\Sp(W)\end{pmatrix}\subset\GL(F\oplus W).\]
It suffices to show that the metaplectic $S^{1}$-cover $\wt{G}$ of $G:=\Aut_{F}(H(W),\sigma)$ reduces to a double cover. Let $[\wt{G}]\in H^{2}(G,S^{1})$ be the cohomology class corresponding to the $S^{1}$-cover $\wt{G}$. Let $i:U\hrar G$ be the inclusion of the unipotent radical of $G$ into $G$, so that $G/U\cong \Sp(W)$. For any choice of a polarization of $H(W)$, it is preserved by $U$, as it only changes the center. Therefore, $i^{*}[\wt{G}]$ is trivial. Moreover, $H^{1}(\Sp(W),S^{1})=\Hom(\Sp(W),S^{1})=\lbrace1\rbrace$ as $\Sp(W)$ is perfect. Therefore, by the Hochschild--Serre spectral sequence (e.g., \cite[Theorem 2]{Hitta}), $[\wt{G}]=q^{*}\alpha$ for some $\alpha\in H^{2}(\Sp(W),S^{1})$, where $q:G\thrar\Sp(W)$ is the quotient map. It is well-known that $\alpha$ is induced from $\alpha'\in H^{2}(\Sp(W),\lbrace\pm1\rbrace)$ corresponding to the metaplectic double cover $\wh{\Sp(W)}\thrar\Sp(W)$ via the inclusion $\lbrace\pm1\rbrace\hrar S^{1}$. Therefore, various functorialities of the smooth cohomology groups imply that $q^{*}\alpha'\in H^{2}(G,\lbrace\pm1\rbrace)$ induces $[\wt{G}]$ via the inclusion $\lbrace\pm1\rbrace\hrar S^{1}$, which implies that the $S^{1}$-cover $\wt{G}$ reduces to the pullback of the metaplectic double cover $q^{*}\wh{\Sp(W)}$, which is a double cover of $G$, as desired.
\end{enumerate}
\end{proof}

We certainly expect that there should be many examples of $\deg\wh{S_{\rho_{\psi}}}>2$ when $\depth\psi\ge3$, although at the moment we are not aware of any.

Here is the simplest example of $\depth\psi=3$, arising from a Borel subgroup of $\GL_{4}(F)$.

\begin{exam}\label{exam:GL4}
Consider $\bfG=\GL_{4}$ over $F$. Let $\bfP\le\bfG$ be the standard upper triangular Borel subgroup, $\bfN$ be its unipotent radical, and $N=\bfN(F)$. We may identify $\Lie N$ with the upper triangular nilpotent matrices. For $a\in F^{\times}$ and $b\in F$, consider $\psi_{a,b}:\Lie N\rar\bC$ defined as
\[\psi_{a,b}\begin{pmatrix}
0 & x_{1,2} & x_{1,3} & x_{1,4}
\\
0 & 0 & x_{2,3} & x_{2,4}
\\
0 & 0 & 0 & x_{3,4}
\\
0 & 0 & 0 & 0
\end{pmatrix}:= ax_{1,4}+bx_{2,3}.\]
As $\psi_{a,b}$ is nontrivial at the top-right corner, it is clear that $\depth\psi_{a,b}=3$. The coadjoint orbit $\cO_{\psi_{a,b}}$ is four-dimensional; if we identify
\[(\Lie N)^{*}=\left\lbrace\begin{pmatrix}
* & * & * & *
\\
y_{1,2} & * & * & *
\\
y_{1,3} & y_{2,3} & * & *
\\
y_{1,4} & y_{2,4} & y_{3,4} & *
\end{pmatrix}\right\rbrace,\]
 \[\left\langle \begin{pmatrix}
* & * & * & *
\\
y_{1,2} & * & * & *
\\
y_{1,3} & y_{2,3} & * & *
\\
y_{1,4} & y_{2,4} & y_{3,4} & *
\end{pmatrix},\begin{pmatrix}
0 & x_{1,2} & x_{1,3} & x_{1,4}
\\
0 & 0 & x_{2,3} & x_{2,4}
\\
0 & 0 & 0 & x_{3,4}
\\
0 & 0 & 0 & 0
\end{pmatrix}\right\rangle=\sum_{1\le i<j\le 4}x_{i,j}y_{i,j},\]then,
\[\cO_{\psi_{a,b}}=
\left\lbrace
\begin{pmatrix}
* & * & * & *
\\
y_{1,2} & * & * & *
\\
y_{1,3} & y_{2,3} & * & *
\\
a & y_{2,4} & y_{3,4} & *
\end{pmatrix}
~:~a(y_{2,3}-b)=y_{1,3}y_{2,4}
\right\rbrace.
\]
Note also that
\[S_{\rho_{\psi_{a,b}}}=\begin{cases}
\lbrace\diag(x,y,y,x)~:~x,y\in F^{\times}\rbrace
&\text{ if }b\ne0,
\\
\lbrace\diag(x,y,z,x)~:~x,y,z\in F^{\times}\rbrace
&\text{ if }b=0.
\end{cases}
\]Note that $S_{\rho_{\psi_{a,b}}}$ stabilizes  a polarization
\[\fh=\left\lbrace\begin{pmatrix}
0 & 0 & 0 & x_{1,4}
\\
0 & 0 & x_{2,3} & x_{2,4}
\\
0 & 0 & 0 & x_{3,4}
\\
0 & 0 & 0 & 0
\end{pmatrix}\right\rbrace\subset\Lie N\]subordinate to $\psi_{a,b}$, which by Lemma~\ref{lem:split over parabolic} implies that $\wh{S_{\rho_{\psi_{a,b}}}}=S_{\rho_{\psi_{a,b}}}$ is of degree $1$.
\end{exam}

We now prove Theorem~\ref{thm:finite cover exists}.
\begin{proof}[Proof of Theorem~\ref{thm:finite cover exists}]

In terms of continuous cohomology, we would like to show that the cohomology class $[\wt{S_{\rho_{\psi}}}]\in H^{2}(S_{\rho_{\psi}},S^{1})$ corresponding to the central extension $\wt{S_{\rho_{\psi}}}$ 
is induced from a cohomology class from $H^{2}(S_{\rho_{\psi}},\mu_{n})$ for some $n\ge1$, where the cohomology groups are the continuous group cohmology groups, and the coefficient modules have the trivial actions. By the short exact sequence $1\rar\mu_{n}\rar S^{1}\xrar{z\mapsto z^{n}} S^{1}\rar1$, it suffices to show that $[\wt{S_{\rho_{\psi}}}]$ is torsion. If we denote $U:=R_{u}(S_{\rho_{\psi}})$ and $R:=S_{\rho_{\psi}}/U$, by the Hochschild--Serre spectral sequence for $U\trianglelefteq S_{\rho_{\psi}}$, we have $H^{i}(R,H^{j}(U,S^{1}))\Rar H^{2}(S_{\rho_{\psi}},S^{1})$. We claim that $H^{2}(S_{\rho_{\psi}},S^{1})\cong H^{2}(R,S^{1})$. For this, we need to show that $H^{2}(U,S^{1})^{R}=H^{1}(R,H^{1}(U,S^{1}))=0$, which both follow from the same argument as in \cite[\S3, Proposition]{Sury}. 

Similarly, there is a Hochschild--Serre spectral sequence $H^{i}(R/Z(R),H^{j}(Z(R),S^{1}))\Rar H^{i+j}(R,S^{1})$. As $H^{2}(Z(R),S^{1})$ is a torsion group by \cite[Corollary 4.1]{Moore}, there exists $N>1$ such that $N[\wt{S_{\rho_{\psi}}}]\in \ker(H^{2}(R,S^{1})\rar H^{2}(Z(R),S^{1}))$. As $R/Z(R)$ is a semisimple $p$-adic group, $\frac{R/Z(R)}{[R/Z(R),R/Z(R)]}$ is finite. Therefore, $H^{1}(R/Z(R),A)=\Hom(R/Z(R),A)$ is a torsion group for any trivial $R/Z(R)$-module $A$. Therefore, there exitss $M>1$ such that $M[\wt{S_{\rho_{\psi}}}]\in\im(H^{2}(R/Z(R),S^{1})\rar H^{2}(R,S^{1}))$. We claim that $H^{2}(R/Z(R),S^{1})$ is a finite group; if this is true, it will finish the proof. Let $R'\thrar R/Z(R)$ be a central isogeny where $R'$ is simply connected. Let $K:=\ker(R'\thrar R/Z(R))$. Then, there is a Hochschild--Serre spectral sequence $H^{i}(R/Z(R),H^{j}(K,S^{1}))\Rar H^{i+j}(R',S^{1})$. By \cite{Moore,PrasadRaghunathan1, PrasadRaghunathan2, PrasadRaghunathan3}, it is known that $H^{2}(R',S^{1})$ is a finite group. As $H^{1}(R',S^{1})=0$, to show that $H^{2}(R/Z(R),S^{1})$ is a finite group, it suffices to show that $H^{1}(K,S^{1})^{R/Z(R)}$ is a finite group, but it is a subgroup of $H^{1}(K,S^{1})=\Hom(K,S^{1})$, which surely is a finite group, as desired.
\end{proof}

\begin{defn}\label{defn:rhomu}
Let $\psi\in(\Lie N)^{*}$, and suppose that a metaplectic cover $\wh{S_{\rho_{\psi}}}$ exists (e.g., when Assumption (\ref{tag:reductive}) holds). Let $\wh{R_{\rho_{\psi}}}:=\wh{S_{\rho_{\psi}}}\ltimes N$ be the group where $\wh{S_{\rho_{\psi}}}$ acts on $N$ via the natural action of $M$ and $\wh{S_{\rho_{\psi}}}\thrar S_{\rho_{\psi}}\hrar M$. Let $(\rho_{\psi},V)$ be the irreducible smooth admissible  representation corresponding to $\psi$ by Theorem~\ref{thm:Kirillov}. Choose an extension of $\rho_{\psi}$ to $\wh{R_{\rho_{\psi}}}$ (i.e., an action of $\wh{R_{\rho_{\psi}}}$ on $V$ such that its restriction to $N$ recovers $\rho_{\psi}$), which will again be denoted as $\rho_{\psi}$.  It is clear that $\rho_{\psi}$ is smooth, admissible and irreducible as an $\wh{R_{\rho_{\psi}}}$-representation.

Let $\sigma$ be an irreducible smooth admissible representation of the metaplectic cover $\wh{S_{\rho_{\psi}}}$, which may also be regarded as an $\wh{R_{\rho_{\psi}}}$-representation by trivially extending over $N$. We say that $\sigma$ accompanies $\rho_{\psi}$ if the $\wh{R_{\rho_{\psi}}}$-representation $\sigma\otimes\rho_{\psi}$ factors through an $R_{\rho_{\psi}}:=S_{\rho_{\psi}}\ltimes N$-representation. It is clear that, for such $\sigma$, $\sigma\otimes\rho_{\psi}$ is an irreducible $R_{\rho_{\psi}}$-representation.
\end{defn}

\begin{rmk}
Note that, in Definition~\ref{defn:rhomu}, $\rho_{\psi}$ is not a unique lift of the natural projective representation of $\wh{R_{\rho_{\psi}}}$ on $V$ that coincides with the given linear action of $N$ on $V$: it is only unique up to a twist by a smooth character of $\wh{S_{\rho_\psi}}$. On the other hand, by considering $\sigma\otimes\rho_{\psi}$, we account for all lifts of the projective representation of $\wh{R_{\rho_{\psi}}}$ on $V$, as $\sigma$ can be twisted by a character of $\wh{S_{\rho_{\psi}}}$.
\end{rmk}
\section{{Generalized Whittaker models beyond $\fsl_{2}$-triples}}\label{sec:GWh}
In this section, we define the \emph{generalized Whittaker models} using smooth irreducible representations of $N$ that do not necessarily arise from $\fsl_{2}$-triples; see Remark~\ref{rem:history} for the comparison of our definition with those appearing in the  literature. In particular, we also define the corresponding generalizations of \emph{degenerate Whittaker models} and \emph{Jacquet functors}, and clarify their fundamental properties. These will be useful in the later part of the paper when we study more refined information regarding the generalized / degenerate Whittaker models.

The  generalization of Jacquet functors is based on the following construction. 
\begin{defn}[Largest isotypic quotient]\label{def:coisotype}
Let $H$ be a td-group, and $(\pi,V)$ and $(\rho,W)$ be smooth  representations of $H$. Suppose that $\rho$ is irreducible. Then, the \emph{largest $\rho$-isotypic quotient} (or \emph{$\rho$-co-isotype}) $\pi_{H,\rho}$ (or $\pi_{\rho}$ if $H$ is clear from the context) is the smallest quotient of $(\pi,V)$ such that every $H$-homomorphism $\pi\rar\rho$ factors through $\pi_{H,\rho}$. Concretely, $\pi_{H,\rho}$ is constructed as the natural $H$-representation structure on \[V_{H,\rho}:=\frac{V}{\bigcap_{f\in\Hom_{H}(\pi,\rho)}\ker f}.\]
It is obvious from the definition that the largest $\rho$-isotypic quotient exists and is unique. It is also clear that $\pi_{H,\rho}$ is a smooth representation of $H$. The universal property can be written as
\[\Hom_{H}(\pi,\rho)=\Hom_{H}(\pi_{H,\rho},\rho).\]
If $\rho$ is the trivial representation, we  omit $\rho$ from the notation, and  call $\pi_{H}$ the \emph{$H$-coinvariant} of $\pi$.
\end{defn}

\begin{rmk}
Another point of view on $V_{H,\rho}$ is as follows. There is an $H$-morphism
$\prod_{f\in\Hom_{H}(\pi,\rho)}f:V\rar \prod_{f\in\Hom_{H}(\pi,\rho)}W$.
Then, $V_{H,\rho}$ is defined as the image of $\prod_{f\in\Hom_{H}(\pi,\rho)}f$. In particular, it is clear that $V_{H,\rho}$ is a smooth $H$-subrepresentation of $\prod_{f\in\Hom_{H}(\pi,\rho)}W$, thus a smooth $H$-subrepresentation of $\left(\prod_{f\in\Hom_{H}(\pi,\rho)}W\right)_{\infty}$. 
\end{rmk}

In general, as noted above, the notion of the largest isotypic quotient involves an awfully big infinite direct product, and it is unclear whether the largest $\rho$-isotypic quotient deserves to be called ``$\rho$-isotypic''. We will see shortly that there is in fact a satisfying definition of $\rho$-isotypic representations when $H$ is a $p$-adic unipotent group.

We now define the generalized Whittaker models as follows.

\begin{defn}[Generalized Whittaker models]\label{def:GWh}
Let $\pi$ be a smooth representation of $G$. Let $\psi\in(\Lie N)^{*}$. Suppose that a metaplectic cover $\wh{S_{\rho_{\psi}}}$ exists, and choose an extension of $\rho_{\psi}$ to $\wh{R_{\rho_{\psi}}}$. Let $\sigma$ be an irreducible smooth admissible representation of the metaplectic cover $\wh{S_{\rho_{\psi}}}$ which accompanies $\rho_{\psi}$. The \emph{generalized Whittaker model} of $\pi$ with respect to $(\sigma,\psi)$ is defined as
\[
\GWh_{\sigma,\psi}(\pi):=\Hom_{G}(\pi,\Ind_{R_{\rho_{\psi}}}^{G}(\sigma\otimes\rho_{\psi}))
=\Hom_{R_{\rho_{\psi}}}(\pi,\sigma\otimes\rho_{\psi}).\]
If  $\sigma$ is the trivial representation, we will often omit the corresponding letters from the notation and just write it as $\GWh_{\psi}(\pi)$.

A generalized Whittaker model with respect to $(\sigma,\psi)$  is called \emph{high-depth} (\emph{low-depth}, respectively) if $\rho_{\psi}$ is so. 
\end{defn}

\begin{lem}\label{lem:GWh}
Retain the notations of Definition~\ref{def:GWh}. Then,
\[\GWh_{\sigma,\psi}(\pi)\cong\Hom_{\wh{R_{\rho_{\psi}}}}(\pi\otimes\rho_{\psi}^{\vee},\sigma)\cong\Hom_{\wh{S_{\rho_{\psi}}}}\left(\left(\pi\otimes\rho_{\psi}^{\vee}\right)_{N},\sigma\right).\]
\end{lem}
\begin{proof}
Note that $\Hom_{R_{\rho_{\psi}}}(\pi,\sigma\otimes\rho_{\psi})=\Hom_{\wh{R_{\rho_{\psi}}}}(\pi,\sigma\otimes\rho_{\psi})$. 
As $\sigma\otimes\rho_{\psi}$ is an irreducible $\wh{R_{\rho_{\psi}}}$-representation, it is admissible, so by Corollary~\ref{coro:contragredient-tensor}, $\left(\sigma\otimes\rho_{\psi}\right)^{\vee}\cong\sigma^{\vee}\otimes\rho_{\psi}^{\vee}$. Therefore, by Corollary~\ref{coro:V1V2V3},
\[\GWh_{\sigma,\psi}(\pi)=\Hom_{\wh{R_{\rho_{\psi}}}}(\pi\otimes\left(\sigma\otimes\rho_{\psi}\right)^{\vee},\bfone)=\Hom_{\wh{R_{\rho_{\psi}}}}(\pi\otimes\sigma^{\vee}\otimes\rho_{\psi}^{\vee},\bfone)=\Hom_{\wh{R_{\rho_{\psi}}}}(\pi\otimes\rho_{\psi}^{\vee},\sigma).\]
The second inequality holds as $N$ acts trivially on $\sigma$.
\end{proof}

\begin{rmk}[Comparison with the generalized Whittaker models in the literature]\label{rem:history}
In the literature, the notion of generalized Whittaker\footnote{Some authors prefer to use ``Gelfand--Graev'' in place of ``Whittaker''. Both words may be used interchangeably.} models is defined and used in slightly varying contexts. We will review them and compare with our notion of generalized Whittaker models in this Remark. It is clear that, by Lemma~\ref{lem:depth}, all prior instances of generalized Whittaker models in the literature, as well as the Bessel and Fourier--Jacobi models of \cite[\S14]{GGP}, are special cases of \emph{low-depth} generalized Whittaker models in the sense of Definition~\ref{def:GWh} (i.e., $\depth\le2$). 

There have been largely two different ways to define the notion of generalized Whittaker models. To distinguish between them (and from our Definition~\ref{def:GWh}), we will call them generalized \emph{A-Whittaker} and \emph{B-Whittaker} models below.
\begin{enumerate}[label=(\Alph*)]
\item One way (e.g., \cite{Kawanaka,Yamashita1,GomezZhu,GomezGourevitchSahi}) is to define them as special cases of \emph{degenerate Whittaker models} (see \cite{MoeglinWaldspurger,GomezGourevitchSahi}). To define degnerate Whittaker models, we choose a semisimple element $h\in\Lie G$ whose adjoint action on $\fg$ has eigenvalues in $\bQ$, and $f\in(\Lie G)^{*}$ such that $h\cdot f=-2f$. Let $U\le G$ ($U^{+}\le G$, respectively) be the unipotent subgroup where $\Lie U\le \Lie G$ is the sum of $h$-eigenspaces of eigenvalues $\ge1$ ($\ge2$, respectively). Note that $f$ naturally defines a character $\chi_{f}$ of $U^{+}$ by $\chi_{f}(g)=\varepsilon(f(\log g))$, and it extends to a slightly bigger group $U'$ which is generated by $U^{+}$ and the stabilizer of $\chi_{f}$ in $U$. If $U\ne U'$, then $\frac{U}{U'\cap\ker\chi_{f}}$ is a Heisenberg group with the center $\frac{U'}{U'\cap\ker\chi_{f}}$, so one can take the Weil representation $\omega_{f}$ with central charcater $\chi_{f}$ and regard it as a smooth irreducible representation of $U$.

Then, for a smooth admissible representation $(\pi, V)$ of $G$, \[\text{Degenerate Whittaker model of $\pi$ w.r.t. $(h,f)$}=\begin{cases}\Hom_{U}(\pi,\chi_{f})&\text{ if $U=U',$}\\\Hom_{U}(\pi,\omega_{f})&\text{ if $U\ne U'$.}\end{cases}\]
In this context, \emph{generalized A-Whittaker models} are degenerate Whittaker models where $(h,f)$ arises from an $\fsl_{2}$-triple; namely, after $G$-equivariantly identifying $\Lie G$ and $(\Lie G)^{*}$, there is an element $e\in \Lie G$ such that $(f,h,e)$ forms an $\fsl_{2}$-triple. Note that, if $(h,f)$ arises from an $\fsl_{2}$-triple, then $U=U'$ if and only if $U=U^{+}$.

A generalized A-Whittaker model is, by construction, ought to be large in size, because the space of degenerate Whittaker model admits a natural group action: if $M\le G$ is the stabilizer group of $h$, then $M$ acts on $U$, so the (possibly central cover of) centralizer of either $\chi_{f}$ or $\omega_{f}$ in $M$ acts on the space of degenerate Whittaker model. Using the notations of Definition~\ref{def:GWh}, this version is missing information about the $\sigma$-part. It is therefore unreasonable to expect the space of generalized A-Whittaker model  to be finite-dimensional unless $\pi$ is ``small.''

\item Another way (e.g., \cite{Yamashita2,GanPeng}\footnote{\cite{Yamashita2} used a rather different name: \emph{reduced generalized Gelfand--Graev representations} (RGGGRs).}), which is a special case of Definition~\ref{def:GWh}, deals with the action of the stabilizer in the Levi, still within the context of $\fsl_{2}$-triples. We start by choosing an $\fsl_{2}$-triple $\gamma=\lbrace e,h,f\rbrace\subset\Lie G$, and let $M\le G$ be the centralizer of $\gamma$. Similarly as above, let $U\le G$ ($U^{+}\le G$, respectively) be the unipotent subgroup where $\Lie U\le\Lie G$ is the sum of $h$-eigenspaces of eigenvalue $\ge1$ ($\ge2$, respectively). After $G$-equivariantly identifying $\Lie G\cong (\Lie G)^{*}$, $f$ defines a character $\chi_{f}$ of $U^{+}$, and if $U\ne U^{+}$, we can construct a Weil representation $\omega_{f}$ with central character $\chi_{f}$, regarded as a smooth irreducible representation of $U$. Note that $\chi_{f}$ extends naturally to a character of $MU^{+}$, but, if $U\ne U^{+}$, $\omega_{f}$ only naturally extends to a representation of the double cover $\wh{M}U$ of $MU$ (compare with  Corollary~\ref{coro:deg}). If $U=U^{+}$, let $\sigma$ be a smooth irreducible representation of $M$. If $U\ne U^{+}$, let $\sigma$ be a genuine smooth irreducible representation of $\wh{M}$, so that $\wh{M}\otimes\omega_{f}$ descends to a representation of $MU$.

Then, for a smooth admissible representation $(\pi, V)$ of $G$,
\[\text{Generalized B-Whittaker model of $\pi$ w.r.t. $\gamma$}=\begin{cases}\Hom_{MU}(\pi,\sigma\otimes\chi_{f})
& \text{ if $U=U^{+}$},
\\\Hom_{MU}(\pi,\sigma\otimes\omega_{f})
& \text{ if $U\ne U^{+}$}.
\end{cases}\]
It is easy to see that Fourier--Jacobi and Bessel models of \cite[\S14]{GGP}, which are now known to be of multiplicity one (see \cite{GGP,AGRS,Sun,SZ,LS}), are special cases of generalized B-Whittaker models. Also, it is clear that this is a special case of generalized Whittaker models in the sense of Definition~\ref{def:GWh}: note that any generalized B-Whittaker model is of depth $\le2$ (i.e. \emph{low-depth}).
\end{enumerate}
\end{rmk}
In view of Remark~\ref{rem:history}, it is also useful to define a generalization of degenerate Whittaker models.
\begin{defn}[Degenerate Whittaker models]\label{def:DWh}
Let $\pi$ be a smooth representation of $G$. Let $\psi\in(\Lie N)^{*}$. The \emph{degenerate Whittaker model} of $\pi$ with respect to $\psi$ is defined as
\[\hspace*{-0.3cm}
\DWh_{\psi}(\pi):=\Hom_{N}(\pi,\rho_{\psi}).\]By Corollary~\ref{coro:V1V2V3}, $\DWh_{\psi}(\pi)\cong\Hom_{N}(\pi\otimes\rho_{\psi}^{\vee},\bfone)$. 

A degenerate Whittaker model with respect to $\psi$  is called \emph{high-depth} (\emph{low-depth}, respectively) if $\rho_{\psi}$ is so. 
\end{defn}
As noted above, if a metaplectic cover $\wh{S_{\rho_{\psi}}}$ exists, then, upon a choice of an extension of $\rho_{\psi}$ to $\wh{R_{\rho_{\psi}}}$, $\DWh_{\psi}(\pi)$ is naturally endowed with the action of  $\wh{S_{\rho_{\psi}}}$: for $f\in \Hom_{N}(\pi,\rho_{\psi})$, $v\in \pi$ and $g\in\wh{S_{\rho_{\psi}}}$,
\[\tag{$\ast$}\label{eq:Homrep}(g\cdot f)(v):=g\cdot(f(g^{-1}\cdot v)).\]It is straightforward to check that $\Hom_{N}(\pi,\rho_{\psi})\cong\Hom_{N}(\pi\otimes\rho_{\psi}^{\vee},\bfone)$ respects the $\wh{S_{\rho_{\psi}}}$-action on both sides, defined using (\ref{eq:Homrep}). One can obtain the generalized Whittaker model from the degenerate Whittaker model as follows.
\begin{lem}\label{lem:GWhDWh}
Retain the notations of Definition~\ref{def:GWh}. Then,
\[\GWh_{\sigma,\psi}(\pi)\cong\Hom_{\wh{S_{\rho_{\psi}}}}(\sigma^{\vee},\DWh_{\psi}(\pi)).\]
\end{lem}
\begin{proof}
This follows from Corollary~\ref{coro:random}, which is applicable as  $\sigma\otimes\rho_{\psi}$ is admissible as an $\wh{R_{\rho_{\psi}}}$-representation. 
\end{proof}

There is a closely related functor $\Rep^{\sm}(G)\rar\Rep^{\sm}(\wh{S_{\rho_{\psi}}})$, which deserves to be called the \emph{Jacquet functor} for generalized Whittaker model. We first observe that the largest $\rho_{\psi}$-isotypic quotient of a smooth $G$-representation has a natural $S_{\rho_{\psi}}\ltimes N$-action (not its metaplectic cover!).
\begin{lem}\label{lem:action of S}
Let $\psi\in(\Lie N)^{*}$. Suppose that a metaplectic cover $\wh{S_{\rho_{\psi}}}$ exists. For a smooth $G$-representation $\pi$, the largest $\rho_{\psi}$-isotypic quotient $\pi_{N,\rho_{\psi}}$ of $\pi$ has a natural $R_{\rho_{\psi}}=S_{\rho_{\psi}}\ltimes N$-action. In other words, the kernel 
\[\ker\left(\pi\thrar\pi_{N,\rho_{\psi}}\right)=\bigcap_{f\in\Hom_{N}(\pi,\rho_{\psi})}\ker f\]is stable under the $R_{\rho_{\psi}}$-action.
\end{lem}
\begin{proof}
Let $W_{\psi}$, $V$ be the underlying vector spaces of the representations $\rho_{\psi}$ and $\pi$, respectively. If $v\in \ker(V\thrar V_{N,\rho_{\psi}})$, or equivalently if there exists $f\in\Hom_{N}(\pi,\rho_{\psi})$ such that $f(v)=0$, then, for any $g\in S_{\rho_{\psi}}$, if we define $f^{g}:V\rar W_{\psi}$ to be the linear map $f^{g}(v)=f(g^{-1}v)$, then $f^{g}(gv)=f(v)=0$, and for any $n\in N$ and $v'\in V$, \[f^{g}(nv')=f(gnv')=f(gng^{-1}\cdot gv')=(gng^{-1})\cdot f(gv')=(gng^{-1})\cdot f^{g}(v').\]
Thus, if we define $\rho_{\psi}^{g}:N\rar\GL(W_{\psi})$ to be a representation $\rho_{\psi}^{g}(n)=\rho_{\psi}(gng^{-1})$, then $f^{g}\in\Hom_{N}(\pi,\rho_{\psi}^{g})$. As $\rho_{\psi}^{g}\cong \rho_{\psi}$ as $N$-representations, this gives that $gv\in\ker(V\thrar V_{N,\rho_{\psi}})$.
\end{proof}
\begin{rmk}\label{rem:not equivariant}
Lemma~\ref{lem:action of S} shows that, if we retain the same notation as in the Lemma, the natural map $\pi_{N,\rho_{\psi}}\hrar\prod_{f\in\Hom_{N}(\pi,\rho_{\psi})}\rho_{\psi}$ is only $N$-equivariant and \emph{not $\wh{R_{\rho_{\psi}}}$-equivariant} (in fact, the action of $R_{\rho_{\psi}}$ on $\pi_{N,\rho_{\psi}}$ does not depend on the choice of an extension of $\rho_{\psi}$ to $\wh{R_{\rho_{\psi}}}$). We defined the $R_{\rho_{\psi}}$-action on $\pi_{N,\rho_{\psi}}$ as above so that the natural isomorphism
\[\Hom_{N}(\pi,\rho_{\psi})\cong\Hom_{N}(\pi_{N,\rho_{\psi}},\rho_{\psi})\]is an isomorphism of $\wh{S_{\rho_{\psi}}}$-representations.
\end{rmk}
\begin{defn}[Jacquet functor for generalized Whittaker model]\label{def:Jacquet}
Let $\psi\in(\Lie N)^{*}$. Suppose that a metaplectic cover $\wh{S_{\rho_{\psi}}}$ exists, and choose an extension of $\rho_{\psi}$ to $\wh{R_{\rho_{\psi}}}$. By Lemma~\ref{lem:action of S}, for a smooth $G$-representation $(V,\pi)$, $\pi_{N,\rho_{\psi}}$  has a natural $R_{\rho_{\psi}}=S_{\rho_{\psi}}\ltimes N$-action (so also an $\wh{R_{\rho_{\psi}}}$-action). Then, the \emph{Jacquet functor} $J_{\rho_\psi}:\Rep^{\sm}(G)\rar\Rep^{\sm}(\wh{S_{\rho_{\psi}}})$ is defined as 
\[J_{\rho_{\psi}}(\pi):=\Hom_{N}(\rho_{\psi},\pi_{N,\rho_{\psi}}),\]
where the action of $\wh{S_{\rho_{\psi}}}$ on $J_{\rho_{\psi}}(\pi)$ is defined using (\ref{eq:Homrep}). 

If we would like to specify the group $N$, we will also use the notation $J_{N,\rho_{\psi}}(\pi)$.
\end{defn}
It is clear that the functor $J_{\rho_{\psi}}$ is well-defined (i.e. $J_{\rho_{\psi}}(\pi)$ is smooth). Also, if $\deg\wh{S_{\rho_{\psi}}}>1$ and $J_{\rho_{\psi}}(\pi)\ne0$, then $J_{\rho_{\psi}}(\pi)$ is always a \emph{genuine} $\wh{S_{\rho_{\psi}}}$-representation.

\begin{rmk}
When $N$ is a Heisenberg group and $\rho_{\psi}$ is a Weil representation, the Jacquet functor as in Definition~\ref{def:Jacquet} as well as its properties proved below was defined and proved in \cite{Weissman} under the name of \emph{Fourier--Jacobi functor}. In \emph{op.~cit.}, it was crucial that a Weil representation is characterized by its central character. This is no longer true for general $p$-adic unipotent groups (see Example~\ref{exam:GL4}), so we obtain analagous results in full generality using deeper results on the smooth representation theory of $p$-adic unipotent groups by \cite{Campbell,Gelfand-Kazhdan}. The richer understanding of smooth representation theory of $p$-adic unipotent group ultimately comes from the fact that any $p$-adic unipotent group is an increasing union of compact open subgroups.
\end{rmk}

We first show that the Jacquet functor for generalized Whittaker model satisfies very similar properties enjoyed by the usual Jacquet functor. We prepare by showing that there is a satisfying definition of what it means to be $\rho$-isotypic for $\rho\in\wh{N}$.

\begin{defn}
Let $\rho\in\wh{N}$. A smooth representation $\pi$ of $N$ is \emph{$\rho$-isotypic} if $\pi=\rho^{\oplus I}$ for some (possibly infinite) index set $I$. Let $\Rep^{\sm}(N,\rho)\subset\Rep^{\sm}(N)$ be the strictly full subcategory of $\rho$-isotypic smooth representations of $N$. By \cite[Proposition 2.10]{Bernstein-Zelevinsky-GLn}, $\Rep^{\sm}(N,\rho)$ is an abelian subcategory of $\Rep^{\sm}(N)$.
\end{defn}

\begin{lem}\label{lem:isotypic}
Let $\pi$ be a smooth representation of $N$. Then, the following are equivalent:
\begin{enumerate}
\item $\pi$ is $\rho$-isotypic;
\item $\pi$ is an $N$-subrepresentation of $\prod_{i\in I}\rho$ (or equivalently $\left(\prod_{i\in I}\rho\right)_{\infty}$) for some (possibly infinite) index set $I$;
\item every irreducible subquotient of $\pi$ is isomorphic to $\rho$.
\end{enumerate}
\end{lem}
\begin{proof}\hfill
\begin{itemize}
\item (1)$\Rightarrow$(2): Suppose that $(\pi,V)$ is $\rho$-isotypic. Then, $\pi\cong\rho^{\oplus I}$ for some index set $I$. It is clear that this is a subrepresentation of $\prod_{i\in I}\rho$.
\item (2)$\Rightarrow$(3): It suffices to show that any irreducible subquotient of $\left(\prod_{i\in I}\rho\right)_{\infty}$ is isomorphic to $\rho$, which is exactly \cite[Lemma 4.2.2]{Campbell}.
\item (3)$\Rightarrow$(1): This follows from \cite[Lemma 2.4.1]{Campbell}.
\end{itemize}
\end{proof}

\begin{coro}\label{coro:tensordecomp}
Let $\psi\in (\Lie N)^{*}$. Suppose that a metaplectic cover $\wh{S_{\rho_{\psi}}}$ exists, and choose an extension of $\rho_{\psi}$ to $\wh{R_{\rho_{\psi}}}$. Let $(\pi,V)$ be a smooth $G$-representation. Then, as $\wh{R_{\rho_{\psi}}}$-representations,
\[\pi_{N,\rho_{\psi}}\cong \rho_{\psi}\otimes J_{\rho_{\psi}}(\pi).\]
\end{coro}
\begin{proof}
Let $W_{\psi}$ be the underlying vector space of the representation $\rho_{\psi}$. There is a natural map $\rho_{\psi}\otimes J_{\rho_{\psi}}(\pi)\rar \pi_{N,\rho_{\psi}}$, $(w,f)\mapsto f(w)$, for $w\in W_{\psi}$ and $f\in J_{\rho_{\psi}}(\pi)=\Hom_{N}(\rho_{\psi},\pi_{N,\rho_{\psi}})$. By \cite[Lemma 2.4.1]{Campbell}, this is an isomorphism of $N$-representations. It is easy to check that this map is also $\wh{S_{\rho_{\psi}}}$-equivariant. Therefore, this map is an isomorphism of $\wh{R_{\rho_{\psi}}}$-representations. 
\end{proof}

\begin{coro}\label{coro:J exact}
Let $\psi\in (\Lie N)^{*}$. Suppose that a metaplectic cover $\wh{S_{\rho_{\psi}}}$ exists, and choose an extension of $\rho_{\psi}$ to $\wh{R_{\rho_{\psi}}}$. 
\begin{enumerate}
\item The largest $\rho_{\psi}$-isotypic quotient functor $(\cdot)_{N,\rho_{\psi}}:\Rep^{\sm}(G)\rar\Rep^{\sm}(N,\rho_{\psi})$ is exact.
\item The Jacquet functor $J_{\rho_{\psi}}:\Rep^{\sm}(G)\rar\Rep^{\sm}(\wh{S_{\rho_{\psi}}})$ is exact.
\end{enumerate}
\end{coro}
\begin{proof}
(2) follows from \cite[Lemma 2.4.2]{Campbell} as the exactness can be checked on the level of vector spaces. As  any $\bC$-vector space is a projective $\bC$-module, (1) then follows from (2) and  Corollary~\ref{coro:tensordecomp}.
\end{proof}

\begin{coro}\label{coro:DWh and J}
Let $\psi\in (\Lie N)^{*}$. Suppose that a metaplectic cover $\wh{S_{\rho_{\psi}}}$ exists, and choose an extension of $\rho_{\psi}$ to $\wh{R_{\rho_{\psi}}}$. For a smooth representation $\pi$ of $G$, we have a natural isomorphism of smooth $\wh{S_{\rho_{\psi}}}$-representations
\[\DWh_{\rho_{\psi}}(\pi)_{\infty}\cong J_{\rho_{\psi}}(\pi)^{\vee}.\]
Therefore,
\[\GWh_{\sigma,\psi}(\pi)\cong\Hom_{\wh{S_{\rho_{\psi}}}}(\sigma^{\vee},J_{\rho_{\psi}}(\pi)^{\vee})\cong\Hom_{\wh{S_{\rho_{\psi}}}}(J_{\rho_{\psi}}(\pi),\sigma),
\]
for any irreducible smooth admissible representation $\sigma$ of $\wh{S_{\rho_{\psi}}}$ accompanying $\rho_{\psi}$.
\end{coro}
\begin{proof}
We claim that, as $\wh{S_{\rho_{\psi}}}$-representations, $\DWh_{\rho_{\psi}}(\pi)\cong\Hom(J_{\rho_{\psi}}(\pi),\bC)$. Note that, by Corollary~\ref{coro:tensordecomp},
\[\DWh_{\rho_{\psi}}(\pi)\cong\Hom_{N}(\pi_{N,\rho_{\psi}},\rho_{\psi})\cong\Hom_{N}(\rho_{\psi}\otimes J_{\rho_{\psi}}(\pi),\rho_{\psi}).\]
There is a natural map $\iota:\Hom(J_{\rho_{\psi}}(\pi),\bC)\rar\Hom_{N}(\rho_{\psi}\otimes J_{\rho_{\psi}}(\pi),\rho_{\psi})$ given by, for $f:J_{\rho_{\psi}}(\pi)\rar\bC$, $v\in\rho_{\psi}$ and $w\in J_{\rho_{\psi}}(\pi)$, $\iota(f)(v\otimes w)=f(w)v$. It is easy to check that this map is $\wh{S_{\rho_{\psi}}}$-equivariant, and we claim that this is a bijection. 

Injectivity is clear. To show the surjectivity, suppose that $h:\rho_{\psi}\otimes J_{\rho_{\psi}}(\pi)\rar\rho_{\psi}$ is $N$-equivariant. For any $w\in J_{\rho_{\psi}}(\pi)$, the map $\rho_{\psi}\rar\rho_{\psi}\otimes J_{\rho_{\psi}}(\pi)$, $v\mapsto v\otimes w$, is $N$-equivariant. Thus, $v\mapsto h(v\otimes w)$ is an element of $\Hom_{N}(\rho_{\psi},\rho_{\psi})$. By Schur's lemma, there exists $f(w)\in\bC$ such that $h(v\otimes w)=f(w)v$ for all $v\in\rho_{\psi}$. Obviously $f$ is $\bC$-linear, so $f$ defines an element of $\Hom(J_{\rho_{\psi}}(\pi),\bC)$, and $h=\iota(f)$, as desired.

The last statement follows from Lemma~\ref{lem:GWhDWh}.
\end{proof}

\begin{rmk}
In the context of degenerate Whittaker models in the usual sense (see Remark~\ref{rem:history}), the archimedean analogue of Corollary~\ref{coro:DWh and J} is shown in \cite[Lemma 4.8]{GomezZhu}. 
It is also clear from the above proof that $\DWh_{\rho_{\psi}}(\pi)$ is smooth if and only if either $\rho_{\psi}$ is finite-dimensional or $J_{\rho_{\psi}}(\pi)$ is finite-dimensional.
\end{rmk}
{
\begin{rmk}
In view of Remark~\ref{rem:not equivariant}, one may alternatively define the $\wh{S_{\rho_{\psi}}}$-action on $\pi_{N,\rho_{\psi}}$ to be the one induced as a subrepresentation of $\prod_{f\in\Hom_{N}(\pi,\rho_{\psi})}\rho_{\psi}$. Let us denote this representation as $(\pi_{N,\rho_{\psi}})^{\mathrm{nm}}$ ($\mathrm{nm}$ stands for ``non-metaplectic''). We can then corresondingly define the non-metaplectic versions of degenerate Whittaker model and Jacquet functor,\[\DWh_{\rho_{\psi}}^{\nm}(\pi):=\Hom_{N}((\pi_{N,\rho_{\psi}})^{\nm},\rho_{\psi}),\quad J_{\rho_{\psi}}^{\nm}(\pi):=\Hom_{N}(\rho_{\psi},(\pi_{N,\rho_{\psi}})^{\nm}).\]
It is easy to check that $\DWh_{\rho_{\psi}}^{\nm}(\pi)$ and $J_{\rho_{\psi}}^{\nm}(\pi)$ are in fact $S_{\rho_{\psi}}$-representations. Therefore, for an irreducible smooth admissible representation $\sigma$ of $S_{\rho_{\psi}}$ (not its metaplectic cover!), one can define the non-metaplectic version of the generalized Whittaker model as
\[\GWh_{\sigma,\psi}^{\nm}(\pi):=\Hom_{S_{\rho_{\psi}}}(\sigma^{\vee},\DWh_{\psi}^{\nm}(\pi)).\]
This seems to be genuinely different from the generalized Whittaker model of Definition~\ref{def:GWh}, even in the case of Fourier--Jacobi models, as there is no apparent way to go between irreducible smooth admissible representations of $S_{\rho_{\psi}}$ and genuine such representations of $\wh{S_{\rho_{\psi}}}$ accompanying a specific representation, unless there is a smooth genuine character of $\wh{S_{\rho_{\psi}}}$. We plan to come back to this in a future work.
\end{rmk}
}

The following, which is an analogue of a fundamental homological property of the usual Whittaker model, is an expected consequence of the exactness of the Jacquet functor:

\begin{prop}\label{prop:Ind injective}
Let $\psi\in(\Lie N)^{*}$. Then, $\Ind_{N}^{G}\rho_{\psi}$ is an injective object in $\Rep^{\sm}(G)$.
\end{prop}
\begin{proof}
It suffices to show that $\Hom_{G}(-,\Ind_{N}^{G}\rho_{\psi})=\Hom_{N}((-)\rvert_{N},\rho_{\psi})=\DWh_{\rho_{\psi}}(-)$ is an exact functor. By the proof of Corollary~\ref{coro:DWh and J}, $\DWh_{\rho_{\psi}}(-)\cong\Hom(J_{\rho_{\psi}}(-),\bC)$. As the Jacquet functor $J_{\rho_{\psi}}(-)$ is exact by Corollary~\ref{coro:J exact}(2), $\DWh_{\rho_{\psi}}(-)$ is exact, as desired.
\end{proof}
\begin{rmk}
It is also expected that $\ind_{N}^{G}\rho_{\psi}$ is a projective object in $\Rep^{\sm}(G)$, although it would require studying the finiteness properties of the said representation. Note that, on the other hand, $\rho_{\psi}$ is in general \emph{not} a projective object in $\Rep^{\sm}(N)$, as illustrated by the example in the introduction of \cite{Gelfand-Kazhdan}. 
\end{rmk}

Here is a very curious corollary.

\begin{coro}\label{coro:curious}
Let $\psi\in(\Lie N)^{*}$. Suppose that a metaplectic cover $\wh{S_{\rho_{\psi}}}$ exists. Suppose further that $\deg\wh{S_{\rho_{\psi}}}>1$. Then, for a smooth representation $\pi$ of $G$, $\dim\DWh_{\rho_{\psi}}(\pi)\ne1$.
\end{coro}
\begin{proof}
Suppose the contrary that $\dim \DWh_{\rho_{\psi}}(\pi)=1$. By Corollary~\ref{coro:DWh and J}, $J_{\rho_{\psi}}(\pi)=\chi$ is a genuine character of $\wh{S_{\rho_{\psi}}}$. By Corollary~\ref{coro:tensordecomp}, $\pi_{N,\rho}\cong\rho_{\psi}\otimes\chi$, which descends to an $S_{\rho_{\psi}}$-representation. This means that $\rho_{\psi}\otimes\chi$ gives an extension of the $N$-action on $\rho_{\psi}$ to an $S_{\rho_{\psi}}$-action on $\rho_{\psi}$, which is impossible as $\deg\wh{S_{\rho_{\psi}}}>1$.
\end{proof}

We now address the effect of the action of $P$ on $\wh{N}$.

{

\begin{lem}\label{lem:conjugate isotypic}
Let $\psi\in(\Lie N)^{*}$. Suppose that a metaplectic cover $\wh{S_{\rho_{\psi}}}$ exists. For $g\in P$, define $\psi^{g}\in(\Lie N)^{*}$ as
\[\psi^{g}(\log n):=\psi(\log(g^{-1}ng)),\quad n\in N.\]
Let $(\pi,V)$ be a smooth representation of $G$. Then, $V\rar V$, $v\mapsto g\cdot v$, gives rise to a bijection $V_{N,\rho_{\psi}}\riso V_{N,\rho_{\psi^{g}}}$.
\end{lem}
\begin{proof}
Let $\fh\subset\Lie N$ be a polarization subordinate to $\psi$, $H=\exp(\fh)$, and $\chi_{\psi}:H\rar S^{1}$ be defined as in Theorem~\ref{thm:Kirillov}. Then, $\rho_{\psi}\cong\Ind_{H}^{N}\chi_{\psi}$. First note that $\Lie(gHg^{-1})$ is a polarization subordinate to $\psi^{g}$, because, for any $h_{1},h_{2}\in H$,
\begin{align*}\psi^{g}([\log(gh_{1}g^{-1}),\log(gh_{2}g^{-1})])&=\psi^{g}(\log([gh_{1}g^{-1},gh_{2}g^{-1}]))
\\&=\psi^{g}(\log(g[h_{1},h_{2}]g^{-1}))=\psi(\log([h_{1},h_{2}]))=0,\end{align*}where the notation $[,]$ is used above for both the Lie bracket in $\Lie N$ and the commutator in $N$. Then, $\chi_{\psi^{g}}:gHg^{-1}\rar S^{1}$ satisfies the identity $\chi_{\psi^{g}}(ghg^{-1})=\chi_{\psi}(h)$ for $h\in H$.

We use the realization of $\rho_{\psi}\cong\Ind_{H}^{N}\chi_{\psi}$ as
\[W_{\psi}:=\lbrace f:N\rar\bC\text{ smooth}\mid f(hn)=\chi_{\psi}(h)f(n)\rbrace,\]and similarly for $\rho_{\psi^{g}}\cong\Ind_{gHg^{-1}}^{N}\chi_{\psi^{g}}$, whose realization we denote by $W_{\psi^{g}}$. Given $f\in W_{\psi^{g}}$, we define $f^{g^{-1}}\in W_{\psi}$ as $f^{g^{-1}}(n):=f(gng^{-1})$. This is well-defined, as, for $h\in H$ and $n\in N$,
\[f^{g^{-1}}(hn)=f(ghng^{-1})=\chi_{\psi^{g}}(ghg^{-1})f(gng^{-1})=\chi_{\psi}(h)f^{g^{-1}}(n).\]
We now claim that the map $V\rar V$, $v\mapsto g\cdot v$, descends to $V_{N,\rho_{\psi}}\rar V_{N,\rho_{\psi^{g}}}$. Suppose that $g\cdot v$ vanishes in the quotient $V_{N,\rho_{\psi^{g}}}$. This means that there exists an $N$-equivariant map $\varphi:V\rar W_{\psi^{g}}$ such that $\varphi(g\cdot v)=0$. We define a linear map $\varphi^{g^{-1}}:V\rar W_{\psi}$ as $v'\mapsto \varphi^{g^{-1}}(v'):=\left(\varphi(g\cdot v')\right)^{g^{-1}}$. As checked above, this is well-defined, i.e. $\varphi^{g^{-1}}(v')\in W_{\psi}$ for $v'\in V$. Also, for $n,n'\in N$ and $v'\in V$,
\begin{align*}\left(\varphi^{g^{-1}}(n\cdot v')\right)(n')&=\left(\varphi(gn\cdot v')\right)(gn'g^{-1})
\\
&=\left((gng^{-1})\cdot\varphi(g\cdot v')\right)(gn'g^{-1})\\
&=\left(\varphi(g\cdot v')\right)(gn'g^{-1}\cdot gng^{-1})\\
&=\left(\varphi(g\cdot v')\right)(gn'ng^{-1})\\
&=\left(\varphi^{g^{-1}}(v')\right)(n'n)\\&=\left(n\cdot \varphi^{g^{-1}}(v')\right)(n').
\end{align*}Thus, $\varphi^{g^{-1}}:V\rar W_{\psi}$ is $N$-equivariant. Obviously $\varphi^{g^{-1}}(v)=0$, so the map $v\mapsto g\cdot v$ indeed descends to $V_{N,\rho_{\psi}}\rar V_{N,\rho_{\psi^{g}}}$. Since the construction $\varphi\mapsto \varphi^{g^{-1}}$ is reversible, the map $V_{N,\rho_{\psi}}\rar V_{N,\rho_{\psi^{g}}}$ is bijective.
\end{proof}

\begin{lem}\label{lem:conjugate J}
Retain the notations of Lemma~\ref{lem:conjugate isotypic} and its proof. Identifying $J_{\rho_{\psi}}(\pi)=\Hom_{N}(W_{\psi},V_{N,\rho_{\psi}})$ and $J_{\rho_{\psi^{g}}}(\pi)=\Hom_{N}(W_{\psi^{g}},V_{N,\rho_{\psi^{g}}})$, there is a bijection $J_{\rho_{\psi}}(\pi)\riso J_{\rho_{\psi^{g}}}(\pi)$ defined by 
\[\varphi:W_{\psi}\rar V_{N,\rho_{\psi}}\mapsto \varphi^{g}:W_{\psi^{g}}\rar V_{N,\rho_{\psi^{g}}},\]
\[\varphi^{g}(f):=g\cdot\left(\varphi\left(f^{g^{-1}}\right)\right).\]
Furthermore, the bijection $J_{\rho_{\psi}}(\pi)\riso J_{\rho_{\psi^{g}}}(\pi)$ is equivariant under the identification $\wh{S_{\rho_{\psi^{g}}}}=g\wh{S_{\rho_{\psi}}}g^{-1}$; namely, for $h\in\wh{S_{\rho_{\psi}}}$ and $\varphi\in J_{\rho_{\psi}}(\pi)$, $(h\cdot \varphi)^{g}=(ghg^{-1})\cdot \left(\varphi^{g}\right)$.
\end{lem}
\begin{proof}
Note that the equation makes sense, as it was shown in the proof of Lemma~\ref{lem:conjugate isotypic} that $f^{g^{-1}}\in W_{\psi}$ for $f\in W_{\psi^{g}}$, and $g\cdot(-)$ sends a vector in $V_{N,\rho_{\psi}}$ to a vector in $V_{N,\rho_{\psi^{g}}}$. 
Note that, for $n,n'\in N$,
\[(n\cdot f)^{g^{-1}}(n')=(n\cdot f)(gn'g^{-1})=f(gn'g^{-1}n)=f^{g^{-1}}(n'g^{-1}ng)=\left((g^{-1}ng)\cdot f^{g^{-1}}\right)(n'),\]
so
\[\varphi^{g}(n\cdot f)=g\cdot \left(\varphi\left((n\cdot f)^{g^{-1}}\right)\right)=g\cdot\left(\varphi\left((g^{-1}ng)\cdot f^{g^{-1}}\right)\right)=(ng)\cdot\left(\varphi\left(f^{g^{-1}}\right)\right)=n\cdot\left(\varphi^{g}(f)\right),\]which implies that $\varphi^{g}$ is $N$-equivariant. It is easy to see that this construction is reversible, and we obtain a bijection $J_{\rho_{\psi}}(\pi)\riso J_{\rho_{\psi^{g}}}(\pi)$. The rest is also easily verified, which is left as an exercise to the reader.
\end{proof}

\section{{High-depth multiplicity one examples}}\label{sec:exam}
In this section, we give some examples of high-depth generalized Whittaker models  of multiplicity one, i.e. examples of $(\sigma,\psi)$ with $\depth\psi\ge3$ such that $\dim\GWh_{\sigma,\psi}(\pi)\le1$. By Remark~\ref{rem:history}, these generalized Whittaker models do not come from the usual notion of generalized Whittaker models. We will also observe that in certain cases, the generalized Whittaker model fails to be finite-dimensional as $S$ is ``not large enough''.

{We first show that, for $G=\GL_{n}(F)$, the generalized Whittaker model associated with a certain class of coadjoint orbits of maximal dimension of the maximal unipotent subgroup of $G$ is of multiplicity one for irreducible supercuspidal representations of $G$. 

\begin{thm}\label{thm:GLn}
Let $G=\GL_{n}(F)$ and $P\le G$ be the standard upper triangular Borel subgroup. Let $N\le P$ be the unipotent radical of $P$. For $a\in F^{\times}$, consider $f\in(\Lie N)^{*}$ defined by
\[f\begin{pmatrix} 0 & x_{1,2} & \cdots & x_{1,n-1} & x_{1,n} \\
0 & 0 & \cdots & x_{2,n-1} & x_{2,n} \\
\vdots & \vdots & \ddots & \vdots & \vdots \\
0 & 0 & \cdots & 0 & x_{n-1,n} \\
0 & 0 & \cdots & 0 & 0\end{pmatrix}:=ax_{1,n}
.\]
Let $\rho_{f}$ be the smooth admissible irreducible  representation of $N$ associated to $f$ via Kirllov's orbit method, Theorem~\ref{thm:Kirillov}. Let $\sigma:S_{\rho_{f}}\rar\bC^{\times}$ be a character, where \[S_{\rho_{f}}=\wh{S_{\rho_{f}}}=\lbrace\diag(x,y_{1},\cdots,y_{n-2},x)~:~x,y_{1},\cdots,y_{n-2}\in F^{\times}\rbrace.\]
Then, for any irreducible supercuspidal representation $\pi$ of $\GL_{n}(F)$,\[\dim\GWh_{\sigma,f}(\pi)=\begin{cases} 1& \text{ if the central character of $\pi$ is $\sigma\rvert_{Z(\GL_{n}(F))}$,}\\0&\text{ otherwise.}\end{cases}\]
\end{thm}

\begin{rmk}
Note that, in Theorem~\ref{thm:GLn}, $\depth f=n-1$, which is bigger than $2$ if $n\ge4$ (e.g., if $n=4$, $f=\psi_{a,0}$ of Example~\ref{exam:GL4}). It is also easy to see that $\cO_{f}$ is of maximal dimension among the coadjoint orbits of $N$.
\end{rmk}

\begin{proof}[Proof of Theorem~\ref{thm:GLn}]
We will prove by induction on $n$, extensively using the theory of Gelfand--Kazhdan and Bernstein--Zelevinsky (see \cite{Bernstein-Zelevinsky-Reductive1, Bernstein-Zelevinsky-Reductive2}). 
To  relieve the notation, let $S:=S_{\rho_{f}}$. A polarization subordinate to $f$ can be taken as
\[\fh=\left\lbrace
\begin{pmatrix} 0 & 0 & 0& \cdots & 0 & 0 &*\\
0 & 0 & * & \cdots & * & * & *\\
0 & 0 & 0 & \cdots & * & * & *\\
\vdots & \vdots & \vdots & \ddots & \vdots & \vdots & \vdots \\
0 & 0 & 0 &\cdots & 0 & * & * \\
0 & 0 & 0 &\cdots & 0 & 0 & * \\
0 & 0 & 0 &\cdots & 0 & 0 & 0\end{pmatrix}
\right\rbrace\subset\Lie N,\]i.e., $\fh$ is the set of upper triangular nilpotent matrices whose top row is zero except the rightmost entry. It is an easy calculation that $S$ is indeed as given in the Theorem statement, and as $S$ stabilizes $\fh$, we have $\wh{S}=S$.

Note that $S=Z(\GL_{n}(F))\times S'$, where $S':=\lbrace\diag(1,y_{1},\cdots,y_{n-2},1)~:~y_{1},\cdots,y_{n-2}\in F^{\times}\rbrace$. It is clear that $\GWh_{\sigma,f}(\pi)=0$ if the central character of $\pi$ does not match that of $\sigma$. On the other hand, if they match, then $\GWh_{\sigma,f}(\pi)=\Hom_{S'N}(\pi,\sigma\otimes\rho_{f})$. By Lemma~\ref{lem:GWh},
\[\Hom_{S' N}(\pi,\sigma\otimes\rho_{f})=\Hom_{S'}(\sigma^{\vee},\Hom_{N}(\pi,\rho_{f})).\]Let $H=\exp\fh\subset N$. By the recipe of Theorem~\ref{thm:Kirillov}, $\rho_{f}=\Ind_{H}^{N}(\chi_{f})$, where $\chi_{f}:H\rar\bC^{\times}$ is defined as $\chi_{f}(h):=\varepsilon(f(\log h))$. Therefore, $\Hom_{N}(\pi,\rho_{f})=\Hom_{H}(\pi,\chi_{f})$.

We split $H$ into two parts,
\[H_{1}=\left\lbrace
\begin{pmatrix} 1 & 0 & \cdots & 0 & * \\
0 & 1 & \cdots & 0 & * \\
\vdots & \vdots & \ddots & \vdots & \vdots \\
0 & 0 & \cdots & 1 & * \\
0 & 0 & \cdots & 0 & 1\end{pmatrix}
\right\rbrace,\quad H_{2}=\left\lbrace
\begin{pmatrix} 1 & 0 & \cdots & 0 & 0 \\
0 & 1 & \cdots & * & 0 \\
\vdots & \vdots & \ddots & \vdots & \vdots \\
0 & 0 & \cdots & 1 & 0 \\
0 & 0 & \cdots & 0 & 1\end{pmatrix}
\right\rbrace.\]Namely, $H_{1}$ is the subgroup of $H$ consisting of matrices which coincide with the identity matrix except the rightmost column, and $H_{2}$ is the centrally embedded subgroup of upper triangular unipotent matrices in $\GL_{n-2}(F)$. 

Note that $f$ is trivial on $H_{2}$. Therefore, as $H_{1}$ is a normal subgroup of $H$,  by Corollary~\ref{coro:DWh and J},
\[\Hom_{H}(\pi,\chi_{f})=  \Hom_{H_{2}}(J_{H_{1},\chi_{f}}(\pi),\bfone).\]
As $\pi$ is supercuspidal and irreducible, by \cite[Thm. 4.4]{Bernstein-Zelevinsky-Reductive1}, the restriction $\pi\rvert_{P_{n}}$ to the mirabolic subgroup 
\[P_{n}=\left\lbrace
\begin{pmatrix} * & * & \cdots & * & * \\
* & * & \cdots & * & * \\
\vdots & \vdots & \ddots & \vdots & \vdots \\
* & * & \cdots & * & * \\
0 & 0 & \cdots & 0 & 1\end{pmatrix}
\right\rbrace\]is irreducible and non-degenerate. Therefore, $J_{H_{1},\chi_{f}}(\pi)$ is an irreducible non-degenerate representation of (conjugate of) a smaller mirabolic subgroup\[P_{n-1}=\left\lbrace
\begin{pmatrix} 1 & 0 & \cdots & 0 & 0 \\
* & * & \cdots & * & 0 \\
\vdots & \vdots & \ddots & \vdots & \vdots \\
* & * & \cdots & * & 0 \\
0 & 0 & \cdots & 0 & 1\end{pmatrix}
\right\rbrace.\]
Note that $\sigma$ can be pulled back to a character of $S'H_{2}$, which we still denote by $\sigma$. Then, by Corollary~\ref{coro:random},
\[\Hom_{S'H_{2}}(J_{H_{1},\chi_{f}}(\pi),\sigma)=\Hom_{S'}(\sigma^{\vee},\Hom_{H_{2}}(J_{H_{1},\chi_{f}}(\pi),\bfone)),\]
 where $\sigma:S'H_{2}\rar \bC^{\times}$ is the character $\sigma\rvert_{S'}$ pulled back along $S'H_{2}\thrar S'$.
It is easy to check that the isomorphisms $\Hom_{N}(\pi,\rho_{f})\cong\Hom_{H}(\pi,\chi_{f})\cong\Hom_{H_{2}}(J_{H_{1},\chi_{f}}(\pi),\bfone)$ are $S'$-equivariant. Therefore, 
 \[\GWh_{\sigma,f}(\pi)=\Hom_{S'H_{2}}(J_{H_{1},\chi_{f}}(\pi),\sigma).\]
Let $U_{n-1}$ be the unipotent radical of $P_{n-1}$, i.e.
\[U_{n-1}=\left\lbrace
\begin{pmatrix} 1 & 0 & \cdots & 0 & 0 \\
* & 1 & \cdots & 0 & 0 \\
\vdots & \vdots & \ddots & \vdots & \vdots \\
* & 0 & \cdots & 1 & 0 \\
0 & 0 & \cdots & 0 & 1\end{pmatrix}
\right\rbrace.\]
Let $\psi:U_{n-1}\rar\bC^{\times}$ be a nontrivial character defined by
\[\psi\begin{pmatrix} 1 & 0 & \cdots & 0 & 0 \\
x_{2,1} & 1 & \cdots & 0 & 0 \\
\vdots & \vdots & \ddots & \vdots & \vdots \\
x_{n-1,1} & 0 & \cdots & 1 & 0 \\
0 & 0 & \cdots & 0 & 1\end{pmatrix}:=\varepsilon(x_{2,1}).\]
Then, by  \emph{loc.~cit.}, $J_{H_{1},\chi_{f}}(\pi)=\ind_{P_{n-2}U_{n-1}}^{P_{n-1}}(\tau\otimes\psi)$, where 
\[P_{n-2}=\left\lbrace
\begin{pmatrix} 1 & 0 & 0& \cdots & 0 & 0 \\
0 & 1 & 0 & \cdots & 0 & 0 \\
0 & * & * & \cdots & * & 0 \\
\vdots & \vdots & \vdots & \ddots & \vdots & \vdots \\
0 & * & * & \cdots & * & 0 \\
0 & 0 & 0 & \cdots & 0 & 1\end{pmatrix}
\right\rbrace,\]
and $\tau$ is an irreducible smooth admissible non-degenerate representation of $P_{n-2}$, extended trivially over $U_{n-1}$. 

To study $J_{H_{1},\chi_{f}}(\pi)\rvert_{S'H_{2}}$, we use Mackey theory. For this, note that there are $n-2$ $S'H_{2}$-orbits in $P_{n-2}U_{n-1}\bs P_{n-1}$, as $P_{n-2}U_{n-1}\bs P_{n-1}/S'H_{2}\cong Q_{n-2}\bs\GL_{n-2}(F)/B_{n-2}$, where $B_{n-2}\le \GL_{n-2}(F)$ is the standard upper triangular Borel subgroup, and $Q_{n-2}\le\GL_{n-2}(F)$ is the block lower triangular parabolic subgroup corresponding to the partition $(1,n-3)$. As $Q_{n-2}\bs \GL_{n-2}(F)\cong \bP^{n-3}(F)$ is the projective $(n-3)$-space, there are $n-2$ $B_{n-2}$-orbits in $\bP^{n-3}(F)$, i.e. 
\[\bP^{n-3}(F)=\bA^{n-3}(F)\cup\cdots\cup\bA^{1}(F)\cup\bA^{0}(F).\]
We can take the representatives to be $g_{0},g_{1},\cdots,g_{n-3}$, where $g_{0}=I_{n}$, and, for $1\le k\le n-3$,
\[g_{k}=\begin{pmatrix}
1 & 0 & 0 & 0 \\
0 & 0 & 1 & 0 \\
0 & I_{k} & 0 & 0 \\
0 & 0 & 0 & I_{n-k-2}
\end{pmatrix},\]which is a block square matrix corresponding to the partition $(1,k,1,n-k-2)$. Note that the orbit corresponding to $g_{k}$ is $\bA^{n-3-k}(F)\subset\bP^{n-3}(F)$.

As there are finitely many orbits, by Mackey theory (i.e., a repeated application of \cite[Proposition 1.8]{Bernstein-Zelevinsky-GLn}), there is an increasing filtration $0=W_{0}\subset W_{1}\subset\cdots\subset W_{n-3}\subset W_{n-2}=J_{H_{1},\chi_{f}}(\pi)\rvert_{S'H_{2}}$ of $S'H_{2}$-representations such that 
\[W_{k+1}/W_{k}\cong\ind_{gP_{n-2}U_{n-1}g^{-1}\cap S'H_{2}}^{S'H_{2}}\left((\tau\otimes\psi)^{g}\right),\]where $(\tau\otimes\psi)^{g}(h):=(\tau\otimes\psi)(g^{-1}hg)$.  To ease the notations, let $V_{k}:=P_{n-2}U_{n-1}\cap g_{k}^{-1} S'H_{2}g_{k}$ for $0\le k\le n-3$: note that
\[V_{k}=\left\lbrace\begin{psmatrix}
1 & 0 & 0 & \cdots & 0 & 0 & 0 & \cdots & 0 & 0 & 0 \\
0 & 1 & 0 & \cdots & 0 & 0 & 0 & \cdots & 0 & 0 & 0 \\
0 & * & * & \cdots & * & * & * & \cdots & * & * & 0 \\
0 & * & 0 & \cdots & * & * & * & \cdots & * & * & 0 \\
\vdots & \vdots & \vdots & \ddots & \vdots &\vdots & \vdots & \ddots & \vdots & \vdots & \vdots \\
0 & * & 0 & \cdots & * & * & * & \cdots & * & * & 0 \\
0 & * & 0 & \cdots & 0 & * & * & \cdots & * & * & 0 \\
0 & 0 & 0 & \cdots & 0 & 0 & * & \cdots & * & * & 0 \\
\vdots & \vdots & \vdots & \ddots & \vdots &\vdots & \vdots & \ddots & \vdots & \vdots & \vdots \\
0 & 0 & 0 & \cdots & 0 & 0 & 0 & \cdots & * & * & 0 \\
0 & 0 & 0 & \cdots & 0 & 0 & 0 & \cdots & 0 & * & 0 \\
0 & 0 & 0 & \cdots & 0 & 0 & 0 & \cdots & 0 & 0 & 1 
\end{psmatrix}
\right\rbrace\leftarrow\left(\text{middle row is the $(k+2)$-th row}\right)
.\]
More precisely, $V_{k}$ is the subgroup of $P_{n-2}$ where the middle $(n-3)\times(n-3)$ square matrix on the block diagonal is upper triangular, and the second column can have unspecified entries between the second row and the $(k+2)$-th row.

Note that, for $0\le k\le n-3$, by Frobenius reciprocity (or its Ext-analogue, e.g. \cite[Proposition A.8]{Casselman}), for $i\ge0$,
\begin{align*}
\Ext^{i}_{S'H_{2}}\left(W_{k+1}/W_{k},\sigma\right)
&=\Ext^{i}_{S'H_{2}}\left(\sigma^{-1},\Ind_{g_{k}V_{k}g_{k}^{-1}}^{S'H_{2}}\left(\left((\tau\otimes\psi)^{g_{k}}\right)^{\vee}\right)\right)
\\
&=\Ext^{i}_{g_{k}V_{k}g_{k}^{-1}}\left(\sigma^{-1},\left((\tau\otimes\psi)^{g_{k}}\right)^{\vee}\right)\\
&=\Ext^{i}_{g_{k}V_{k}g_{k}^{-1}}\left((\tau\otimes\psi)^{g_{k}},\sigma\right)\\
&=\Ext^{i}_{V_{k}}\left(\tau\otimes\psi,\sigma^{g_{k}^{-1}}\right)\\
&=\Ext^{i}_{V_{k}}\left(\tau,\sigma^{g_{k}^{-1}}\right),
\end{align*}where the last inequality holds as $V_{k}\cap U_{n-1}=\lbrace1\rbrace$. Note also that $\sigma^{g_{k}^{-1}}$ is still a character which only depends on the diagonal entries of the matrices.

As $\tau$  is an irreducible smooth admissible non-degenerate representation of $P_{n-2}$, by \cite[3.2]{Bernstein-Zelevinsky-Reductive1}, $\tau=\ind_{P_{n-3}U_{n-2}}^{P_{n-2}}(\tau'\otimes\psi')$, where $U_{n-2}\le P_{n-2}$ is the unipotent radical, i.e. 
\[U_{n-2}=\left\lbrace
\begin{pmatrix} 1 & 0 & 0 &\cdots & 0 & 0 \\
0 & 1 & 0 &\cdots & 0 & 0 \\
0 & * & 1 &\cdots & 0 & 0 \\
\vdots & \vdots & \vdots & \ddots & \vdots & \vdots \\
0 & * & 0 &\cdots & 1 & 0 \\
0 & 0 & 0 &\cdots & 0 & 1\end{pmatrix}
\right\rbrace,\]$P_{n-3}\le P_{n-2}$ is the mirabolic subgroup
\[P_{n-3}:=\left\lbrace
\begin{pmatrix} 1 & 0 & 0 &\cdots & 0 & 0 \\
0 & 1 & 0 &\cdots & 0 & 0 \\
0 & 0 & 1 &\cdots & 0 & 0 \\
0 & 0 & * &\cdots & * & 0 \\
\vdots & \vdots & \vdots & \ddots & \vdots & \vdots \\
0 & 0 & * &\cdots & * & 0 \\
0 & 0 & 0 &\cdots & 0 & 1\end{pmatrix}
\right\rbrace,\]
$\tau'$ is an irreducible smooth admissible non-degenerate representation of $P_{n-3}$, extended trivially over $U_{n-2}$, and $\psi':U_{n-2}\rar\bC^{\times}$ is a nontrivial character defined by 
\[\psi'\begin{pmatrix} 1 & 0 & 0 &\cdots & 0 & 0 \\
0 & 1 & 0 &\cdots & 0 & 0 \\
0 & x_{3,2} & 1 &\cdots & 0 & 0 \\
\vdots & \vdots & \vdots & \ddots & \vdots & \vdots \\
0 & x_{n-1,2} & 0 &\cdots & 1 & 0 \\
0 & 0 & 0 &\cdots & 0 & 1\end{pmatrix}:=\varepsilon(x_{3,2}).\]
Again, we would like to use Mackey theory to study $\tau\rvert_{V_{k}}$. Note that, for $V_{0}$, then this situation is exactly analogous to the previous application of Mackey theory, except that the dimension is one lower. Namely, $P_{n-3}U_{n-2}\bs P_{n-2}/V_{0}$ has $n-3$ cosets, whose representatives $h_{0},\cdots,h_{n-4}$ can be taken as $h_{0}=I_{n}$ and, for $1\le l\le n-4$,
\[h_{l}=\begin{pmatrix}
1 & 0 & 0 & 0 & 0 \\
0 & 1 & 0 & 0 & 0 \\
0 & 0 & 0 & 1 & 0 \\
0 & 0 & I_{l} & 0 & 0 \\
0 & 0 & 0 & 0 & I_{n-l-3}
\end{pmatrix},\]which is a block square matrix corresponding to the partition $(1,1,l,1,n-l-3)$. On the other hand, even if $k>0$, as $V_{k}\supset V_{0}$, $h_{0},\cdots,h_{n-4}$ represent all cosets of $P_{n-3}U_{n-2}\bs P_{n-2}/V_{k}$ for all $0\le k\le n-3$ (possibly with duplicates). Therefore, by Mackey theory, $\tau\rvert_{V_{k}}$ has a filtration of $V_{k}$-representations where successive subquotients are of the form $\ind_{h_{l}P_{n-3}U_{n-2}h_{l}^{-1}\cap V_{k}}^{V_{k}}((\tau'\otimes\psi')^{h_{l}})$ for some $0\le l\le n-4$. 

If $k>0$, then $P_{n-3}U_{n-2}\cap h_{l}^{-1}V_{k}h_{l}$ has a nontrivial intersection with $U_{n-2}$, and in particular always has the nontrivial ``$x_{3,2}$-entry''; more precisely, if we let $A_{3,2}\le U_{n-2}$ be the group of lower triangular unipotent matrices with only nonzero off-diagonal entry at the third row and the second column, then $A_{3,2}\subset P_{n-3}U_{n-2}\cap h_{l}^{-1}V_{k}h_{l}$. Therefore, $\psi'\rvert_{\left(P_{n-3}U_{n-2}\cap h_{l}^{-1}V_{k}h_{l}\right)\cap U_{n-2}}$ is nontrivial. On the other hand,  $\left(\sigma^{g_{k}^{-1}}\right)^{h_{l}^{-1}}\rvert_{\left(P_{n-3}U_{n-2}\cap h_{l}^{-1}V_{k}h_{l}\right)\cap U_{n-2}}$ is trivial, as $\left(\sigma^{g_{k}^{-1}}\right)^{h_{l}^{-1}}$ only depends on the diagonal entries. Also, $\tau\rvert_{\left(P_{n-3}U_{n-2}\cap h_{l}^{-1}V_{k}h_{l}\right)\cap U_{n-2}}$ is trivial. Since $\left(P_{n-3}U_{n-2}\cap h_{l}^{-1}V_{k}h_{l}\right)\cap U_{n-2}$ is a normal subgroup of $P_{n-3}U_{n-2}\cap h_{l}^{-1}V_{k}h_{l}$, by the Hochschild--Serre spectral sequence (\cite[Proposition A.9]{Casselman}), we have 
\[\Ext^{i}_{\frac{P_{n-3}U_{n-2}\cap h_{l}^{-1}V_{k}h_{l}}{(P_{n-3}U_{n-2}\cap h_{l}^{-1}V_{k}h_{l})\cap U_{n-2}}}\left(H_{j}((P_{n-3}U_{n-2}\cap h_{l}^{-1}V_{k}h_{l})\cap U_{n-2},\tau'\otimes\psi'),\left(\sigma^{g_{k}^{-1}}\right)^{h_{l}^{-1}}\right)\]\[\Rightarrow\Ext^{i+j}_{P_{n-3}U_{n-2}\cap h_{l}^{-1}V_{k}h_{l}}\left(\tau'\otimes\psi',\left(\sigma^{g_{k}^{-1}}\right)^{h_{l}^{-1}}\right).\]On the other hand, as $(P_{n-3}U_{n-2}\cap h_{l}^{-1}V_{k}h_{l})\cap U_{n-2}$ is an increasing union of compact subsets, by \cite[3.3]{Bernstein-Notes}, $H_{j}((P_{n-3}U_{n-2}\cap h_{l}^{-1}V_{k}h_{l})\cap U_{n-2},-)=0$ for $j>0$. Moreover, $H_{0}((P_{n-3}U_{n-2}\cap h_{l}^{-1}V_{k}h_{l})\cap U_{n-2},\tau'\otimes\psi')=0$ as $(P_{n-3}U_{n-2}\cap h_{l}^{-1}V_{k}h_{l})\cap U_{n-2}$ acts as a nontrivial character on $\tau'\otimes\psi'$ as observed above. Therefore, for $i\ge0$,
\[\Ext^{i}_{V_{k}}\left(\ind_{h_{l}P_{n-3}U_{n-2}h_{l}^{-1}\cap V_{k}}^{V_{k}}((\tau'\otimes\psi')^{h_{l}}),\sigma^{g_{k}^{-1}}\right)=\Ext^{i}_{P_{n-3}U_{n-2}\cap h_{l}^{-1}V_{k}h_{l}}\left(\tau'\otimes\psi',\left(\sigma^{g_{k}^{-1}}\right)^{h_{l}^{-1}}\right)=0.\]
By the Ext long exact sequence, this implies that $\Ext^{i}_{V_{k}}\left(\tau,\sigma^{g_{k}^{-1}}\right)=0$ for all $i\ge0$. This holds for all $k>0$, so again by applying the Ext long exact sequence for the filtration $W_{\bullet}$ of $J_{H_{1},\chi_{f}}(\pi)\rvert_{S'H_{2}}$, we have, for all $i\ge0$, \[\Ext^{i}_{S'H_{2}}(J_{H_{1},\chi_{f}}(\pi),\sigma)=\Ext^{i}_{V_{0}}(\tau,\sigma).\]In particular, $\Hom_{S'H_{2}}(J_{H_{1},\chi_{f}}(\pi),\sigma)=\Hom_{V_{0}}(\tau,\sigma)$. Note that the setting of $\Hom_{V_{0}}(\tau,\sigma)$ is completely analogous to that of $\Hom_{S'H_{2}}(J_{H_{1},\chi_{f}}(\pi),\sigma)$, except that the dimensions of all the involved groups  have been reduced by $1$. Thus, we can inductively apply the same argument as above to show that the above Hom space is equal to $\Hom_{\lbrace1\rbrace}(\bfone,\bfone)$, which is indeed one-dimensional, as desired.
\end{proof}
}
\begin{rmk}\label{rem:GL4}
A similar argument shows that, when $n=4$, $\dim\GWh_{\sigma,\psi_{a,b}}(\pi)=\infty$ when $b\ne0$ and  $\pi$ is irreducible supercuspidal (see Example~\ref{exam:GL4} for the notation). In fact, this is because $\DWh_{\psi_{a,b}}(\pi)$, when $b\ne0$, has a large chunk that is shared with $\DWh_{\psi_{a,0}}(\pi)$, which is naturally a representation of a larger group. See also \S\ref{sec:GL4}.
\end{rmk}
There are also examples of high-depth generalized Whittaker models of multiplicity one for every symplectic group $\Sp_{2n}(F)$ with $n\ge3$ (in $\Sp_{2}(F)=\SL_{2}(F)$ or $\Sp_{4}(F)$, the depth of any unipotent subgroup is $\le2$). 

\begin{defn}\label{def:Sp6}
For $n\ge3$, let $G=\Sp_{2n}(F)$. More precisely, let $D$ be the $n\times n$ anti-diagonal matrix, and $J=\begin{pmatrix}0&D\\-D&0\end{pmatrix}$, and $G=\Sp_{2n}(F)=\left\lbrace A\in\GL_{2n}(F)~:~A^{T}JA=J\right\rbrace$. Let $P\le G$ be the standard upper triangular Borel, i.e., $P=\lbrace A\in \Sp_{2n}(F)~:~A\text{ is upper triangular}\rbrace$. Concretely, 
\[P=\left\lbrace\begin{pmatrix}A& ADB \\ 0& D(A^{T})^{-1}D\end{pmatrix}~:~A,B\in\GL_{n}(F),~\text{$A$ is upper triangular, $B$ is symmetric}\right\rbrace.\]
The Borel subgroup has a Levi decomposition $P=MN$ where \[M=\lbrace\diag(x_{1},\cdots,x_{n},x_{n}^{-1},\cdots,x_{1}^{-1})~:~x_{1},\cdots,x_{n}\in F^{\times}\rbrace,\]
\[N=\left\lbrace\begin{pmatrix}A& ADB \\ 0& D(A^{T})^{-1}D\end{pmatrix}~:~A,B\in\GL_{n}(F),~\text{$A$ is upper triangular unipotent, $B$ is symmetric}\right\rbrace.\]Note that $N$ is a three-step unipotent group. Concretely, $N=U\ltimes N_{Q}$, where $N_{Q}\trianglelefteq N$ is
\[N_{Q}=\left\lbrace\begin{pmatrix}I_{n}& DB \\ 0& I_{n}\end{pmatrix}~:~B\in\GL_{n}(F),~\text{$B$ is symmetric}\right\rbrace,\]and $U\le\GL_{n}(F)$ is the group of upper triangular unipotent $n\times n$ matrices. Note that $N_{Q}$ is the unipotent radical of the \emph{Siegel parabolic} $Q\le G$ defined by
\[Q=\left\lbrace\begin{pmatrix}A& ADB \\ 0& D(A^{T})^{-1}D\end{pmatrix}~:~A,B\in\GL_{n}(F),~\text{$B$ is symmetric}\right\rbrace.\]
We identify the Levi part $M_{Q}$ of $Q$ with $\GL_{n}(F)$ by sending $\begin{pmatrix}A&0\\0&D(A^{T})^{-1}D\end{pmatrix}\mapsto A$.
The Lie algebra $\Lie N$ can be identified with
\[\Lie N=\lbrace A\in M_{n}(F)~:~A^{T}J+JA=0,~\text{$A$ is upper triangular nilpotent}\rbrace.\]
Its center $\fz\subset\Lie N$ is one-dimensional, generated by the top-right corner.
\end{defn}
\begin{thm}\label{thm:Sp6}
Retain the notations of Definition~\ref{def:Sp6}. For $a\in F^{\times}$, consider $f\in(\Lie N)^{*}$ defined by 
\[f\begin{pmatrix} 0 & x_{1,2} & \cdots & x_{1,2n-1} & x_{1,2n} \\
0 & 0 & \cdots & x_{2,2n-1} & x_{2,2n} \\
\vdots & \vdots & \ddots & \vdots & \vdots \\
0 & 0 & \cdots & 0 & x_{2n-1,2n} \\
0 & 0 & \cdots & 0 & 0\end{pmatrix}:=ax_{1,2n}
.\]
Let $\rho_{f}$ be the smooth admissible irreducible  representation of $N$ associated to $f$ via Kirillov's orbit method, Theorem~\ref{thm:Kirillov}. Let $\sigma:S_{\rho_{f}}\rar S^{1}$ be a unitary character, where
\[S_{\rho_{f}}=\wh{S_{\rho_{f}}}=\lbrace\diag(\varepsilon,x_{2},\cdots,x_{n},x_{n}^{-1},\cdots,x_{2}^{-1},\varepsilon)~:~\varepsilon\in\lbrace\pm1\rbrace,~x_{2},\cdots,x_{n}\in F^{\times}\rbrace.\]
Then, for any irreducible smooth admissible representation $\pi$ of $\Sp_{2n}(F)$,
\[\dim\GWh_{\sigma,f}(\pi)\le1.\]
\end{thm}
\begin{rmk}
Note that, in Theorem~\ref{thm:Sp6}, $\depth f=n$. It is also easy to see that $\cO_{f}$ is of maximal dimension among the coadjoint orbits of $N$. On the other hand, the generalized Whittaker model in concern is related to a Fourier--Jacobi model as in \cite[\S14]{GGP} (more specifically, a $(2n-2)$-dimensional symplectic space inside a $2n$-dimensional symplectic space). When $n=2$, the generalized Whittaker model in concern is a Bessel model for $\Sp(4)$ in the sense of \cite{Prasad-Bessel}.
\end{rmk}

\begin{proof}[Proof of Theorem~\ref{thm:Sp6}]
To relieve the notation, let $S:=S_{\rho_{f}}$. It is a straightforward calculation that a polarization subordinate to $f$ can be taken as
\[\fh=\left\lbrace\begin{psmatrix}
0 & 0 & 0 & 0 & \cdots & 0 & 0 & * & * & \cdots & * & * & * & * \\
0 & 0 & * & * & \cdots & * & * & * & * & \cdots & * & * & * & * \\
0 & 0 & 0 & * & \cdots & * & * & * & * & \cdots & * & * & * & * \\
0 & 0 & 0 & 0 & \cdots & * & * & * & * & \cdots & * & * & * & * \\
\vdots & \vdots & \vdots & \vdots & \ddots & \vdots & \vdots & \vdots & \vdots & \ddots & \vdots & \vdots & \vdots & \vdots \\
0 & 0 & 0 & 0 & \cdots & 0 & * & * & * & \cdots & * & * & * & * \\
0 & 0 & 0 & 0 & \cdots & 0 & 0 & * & * & \cdots & * & * & * & * \\
0 & 0 & 0 & 0 & \cdots & 0 & 0 & 0 & * & \cdots & * & * & * & 0 \\
0 & 0 & 0 & 0 & \cdots & 0 & 0 & 0 & 0 & \cdots & * & * & * & 0 \\
\vdots & \vdots & \vdots & \vdots & \ddots & \vdots & \vdots & \vdots & \vdots & \ddots & \vdots & \vdots & \vdots & \vdots \\
0 & 0 & 0 & 0 & \cdots & 0 & 0 & 0 & 0 & \cdots & 0 & * & * & 0 \\
0 & 0 & 0 & 0 & \cdots & 0 & 0 & 0 & 0 & \cdots & 0 & 0 & * & 0 \\
0 & 0 & 0 & 0 & \cdots & 0 & 0 & 0 & 0 & \cdots & 0 & 0 & 0 & 0 \\
0 & 0 & 0 & 0 & \cdots & 0 & 0 & 0 & 0 & \cdots & 0 & 0 & 0 & 0 \\
\end{psmatrix}
\right\rbrace\subset\Lie N.\]
Namely, $\fh$ consists of upper triangular nilpotent matrices in $\Lie N$ such that the first row of the upper-left $n\times n$ block and the last column of the lower-right $n\times n$ block are zero.

It is easy to see that $S$ is indeed as given in the Theorem statement, and as $\fh$ is invariant under the $S$-action, by Lemma~\ref{lem:split over parabolic}, we have $\wh{S}=S$. Let $H=\exp(\fh)\subset N$ and $\chi_{f}:H\rar\bC^{\times}$ is defined as $\chi_{f}(h):=\varepsilon(f(\log h))$. Then, by the recipe of Theorem~\ref{thm:Kirillov}, $\rho_{f}=\Ind_{H}^{N}(\chi_{f})$.
Note that the center of $G$ is $\lbrace\pm I_{2n}\rbrace$, and $S=S'\times\lbrace\pm I_{2n}\rbrace$, where $S'\trianglelefteq S$ is the index $2$ subgroup
\[S':=\lbrace\diag(1,x_{2},\cdots,x_{n},x_{n}^{-1},\cdots,x_{2}^{-1},1)~:~x_{2},\cdots,x_{n}\in F^{\times}\rbrace.\]As $\Hom_{SN}(\pi,\sigma\otimes\rho_{f})\subset\Hom_{S'N}(\pi,\sigma\otimes\rho_{f})$, it suffices to prove that $\dim\Hom_{S'N}(\pi,\sigma\otimes\rho_{f})\le1$.

Let $C\trianglelefteq N$ be the subgroup
\[C:=\left\lbrace\begin{pmatrix}1& * & * & \cdots &  * & * & * \\ 
0 & 1 & 0 & \cdots & 0 & 0 & * \\
0 & 0 & 1 & \cdots & 0 & 0 & * \\
\vdots & \vdots & \vdots & \ddots & \vdots & \vdots & \vdots \\
0 & 0 & 0 & \cdots & 1 & 0 & * \\
0 & 0 & 0 & \cdots & 0 & 1 & * \\
0 & 0 & 0 & \cdots & 0 & 0 & 1\end{pmatrix}\right\rbrace.\]
This is isomorphic to a Heisenberg group. To be more precise, let $V$ be a $2n$-dimensional vector space over $F$ with basis vectors $v_{1},\cdots,v_{2n}$, and let $\langle,\rangle:V\times V\rar F$ be the symplectic form defined by 
\[\left\langle\sum_{i=1}^{2n}a_{i}v_{i},\sum_{i=1}^{2n}b_{i}v_{i}\right\rangle:=\sum_{i=1}^{n} (a_{i}b_{2n+1-i}-a_{n+i}b_{n+1-i}).
\]
Then, $G$ and $\Sp(V)$ are identified in the usual way. Consider the vector subspace $V'\subset V$ spanned by $v_{2},\cdots,v_{2n-1}$, which is endowed with a symplectic form $\langle,\rangle\rvert_{V'}$. Then, $C\cong H(V')$. As $N=H\cdot C$, by Mackey theory, $\rho_{f}\rvert_{C}\cong\Ind_{H\cap C}^{C}(\chi_{f}\rvert_{H\cap C})$. Note that $H\cap C\subset C$ is precisely a polarization corresponding to the Lagrangian subspace $X':=F v_{n+1}+\cdots+F v_{2n-1}\subset V'$. Therefore, $\Ind_{H\cap C}^{C}(\chi_{f}\rvert_{H\cap C})\cong\omega_{\chi_{f}}$, the Weil representation of $C\cong H(V')$ with the central character $\chi_{f}\rvert_{Z(C)}$. If we define
\[N':=\left\lbrace\begin{pmatrix}1& 0 & 0 & \cdots & 0 & 0 & 0 \\ 
0 & 1 & * & \cdots & * & * & 0 \\
0 & 0 & 1 & \cdots & * & * & 0 \\
\vdots & \vdots & \vdots & \ddots & \vdots & \vdots & \vdots \\
0 & 0 & 0 & \cdots & 1 & * & 0 \\
0 & 0 & 0 & \cdots & 0 & 1 & 0 \\
0 & 0 & 0 & \cdots & 0 & 0 & 1\end{pmatrix}\right\rbrace,\]then $S'N'$ stabilizes both $H\cap C$ and $\chi_{f}\cap H\cap C$. Thus, by Lemma~\ref{lem:split over parabolic}, $\omega_{\chi_{f}}$ extends naturally over $S'N$, which is the extension $\rho_{f}$ in the definition of generalized Whittaker model. 

Let $\wh{\Sp(V')}$ be the metaplectic double cover of $\Sp(V')$, and for any subgroup $T\le \Sp(V')$, let $\wh{T}$ be the preimage of $T$ in $\wh{\Sp(V')}$. Let $\omega_{\chi_{f}}^{\wedge}$ be the Weil representation extended as a representation of $\wh{\Sp(V')}\ltimes H(V')$. Let $B_{V'}\le\Sp(V')$ be the upper triangular Borel subgroup (with respect to the choice of basis $\lbrace v_{2},\cdots,v_{2n-1}\rbrace$). Then, we can obtain two extensions of the Weil representation $\omega_{\chi_{f}}$ (a priori a representation of $H(V')\cong C$) as representations of $\wh{B_{V'}}\ltimes H(V')$.
\begin{enumerate}
\item As $\rho_{f}$ is a representation of $S'N=B_{V'}\ltimes H(V')$, so we can pull this representation back to a representation of $\wh{B_{V'}}\ltimes H(V')$. Let us denote this representation  as $\rho_{f}$ by abuse of notation.
\item We can on the other hand restrict the usual Weil representation $\omega_{\chi_{f}}^{\wedge}$ (as a representation of $\wh{\Sp(V')}\ltimes H(V')$) to a subgroup $\wh{B_{V'}}\ltimes H(V')\le \wh{\Sp(V')}\ltimes H(V')$. Let us denote this representation  as $\omega_{\chi_{f}}^{\wedge}$ by abuse of notation.
\end{enumerate}
Both $\rho_{f}$ and $\omega_{\chi_{f}}^{\wedge}$ extend $\omega_{\chi_{f}}$ (i.e., when restricted to $H(V')$, they give $\omega_{\chi_{f}}$), but as representations of $\wh{B_{V'}}\ltimes H(V')$, they are different (for example, $\omega_{\chi_{f}}^{\wedge}$ is genuine, while $\rho_{f}$ is not). On the other hand, by comparing how $\wh{B_{V'}}$ acts using \cite[Lemma 3.2(1)]{Rao}, we know that 
\[\omega_{\chi_{f}}^{\wedge}=\rho_{f}\otimes\chi_{\varepsilon_{a}},\]
where $\chi_{\varepsilon_{a}}(g,\varepsilon):=\varepsilon\gamma(\det g,\frac{1}{2}\varepsilon)^{-1}$; here, $\gamma(a,\eta)$ is the Weil invariant (see \cite[p. 367]{Rao}), and $\wh{B_{V'}}$ is identified, as a set, with $B_{V'}\times\lbrace\pm1\rbrace$, with the multiplication given by $(g_{1},\varepsilon_{1})\cdot(g_{2},\varepsilon_{2})=(g_{1}g_{2},\varepsilon_{1}\varepsilon_{2}(\det g_{1},\det g_{2})_{F})$, where $(-,-)_{F}$ is the quadratic Hilbert symbol in $F$. Therefore,
\begin{align*}\Hom_{S'N}(\pi,\rho_{f}\otimes\sigma)&=\Hom_{B_{V'}\ltimes H(V')}(\pi,\rho_{f}\otimes\sigma)\\&=\Hom_{\wh{B_{V'}}\ltimes H(V')}(\pi,\omega_{\chi_{f}}^{\wedge}\otimes\chi_{\varepsilon_{a}}^{-1}\otimes\sigma)\\&=\Hom_{\wh{\Sp(V')}\ltimes H(V')}\left(\pi,\omega_{\chi_{f}}^{\wedge}\otimes\Ind_{\wh{B_{V'}}\ltimes H(V')}^{\wh{\Sp(V')}\ltimes H(V')}\left(\chi_{\varepsilon_{a}}^{-1}\otimes\sigma\right)\right)\\&=\Hom_{\wh{\Sp(V')}\ltimes H(V')}\left(\pi,\omega_{\chi_{f}}^{\wedge}\otimes\Ind_{\wh{B_{V'}}}^{\wh{\Sp(V')}}\left(\chi_{\varepsilon_{a}}^{-1}\otimes\sigma\right)\right).\end{align*}By \cite[Theorem 3.1]{HanzerMatic}, $\wh{\tau}:=\Ind_{\wh{B_{V'}}}^{\wh{\Sp(V')}}\left(\chi_{\varepsilon_{a}}^{-1}\otimes\sigma\right)$ is a genuine smooth irreducible representation of $\wh{\Sp(V')}$ (the cited Theorem applies, as $\chi_{\varepsilon_{a}}^{-1}=\chi_{\varepsilon_{a}}\otimes\chi_{\varepsilon_{a}}^{-2}$, and $\chi_{\varepsilon_{a}}^{-2}$ is a quadratic character of $B_{V'}$). Therefore, by \cite[Corollary 16.2]{GGP}, $\dim\Hom_{S'N}(\pi,\rho_{f}\otimes\sigma)\le1$, as desired.
\end{proof}

\section{Horizontal degeneration of coadjoint orbits}\label{sec:GL4}
Retain the notations of Definition~\ref{def:GWh}. Lemma~\ref{lem:conjugate J} shows that the Jacquet functors (or, by Corollary~\ref{coro:DWh and J}, the degenerate Whittaker models) for $\psi\in\wh{N}$ in an $M$-orbit contain the same amount of information. Put differently, this means that the nature of the degenerate Whittaker model for a given $\psi\in\wh{N}$ depends on its $M$-orbit $M\psi\subset\wh{N}$. This is a new kind of information not seen on the level of representation theory of $N$, as the ``size'' of the representation $\rho_{\psi}$ is purely determined by the coadjoint $N$-orbit $\cO_{\psi}\subset(\Lie N)^{*}$.

As observed in Section \ref{sec:exam}, The generalized Whittaker model can have vastly different sizes for $N$-coadjoint orbits of equal dimensions. Namely,  there are examples where $\GWh_{\sigma,\psi}(\pi)$ is finite-dimensional while $\GWh_{\sigma,\psi'}(\pi)$ is not, even though $\dim\cO_{\psi}=\dim\cO_{\psi'}$. 
We expect that there is a relation between $\DWh_{\psi}(\pi)$ and $\DWh_{\psi_{0}}(\pi)$ for $\psi,\psi_{0}\in\wh{N}$, when $M\psi_{0}\subset\ov{M\psi}$, so that  $\DWh_{\psi}(\pi)$ is ``large'' partly because it shares a part coming from the degenerate Whittaker model of $\psi_{0}$, which admits an action by a larger stabilizer group.  

We specify the types of degenerations of $M$-orbits we expect to have relationship between degenerate Whittaker models. For this section only, we assume for simplicity that $\bfG$ is split over $F$. There is a split maximal $F$-torus $\bfT\subset\bfM\subset\bfG$ and a parabolic subset $\Psi$ of the set of roots $\Phi(\bfG,\bfT)$ (i.e., $\Phi(\bfG,\bfT)=\Psi\cup-\Psi$, and $\Psi$ is closed under addition) corresponding to the parabolic subgroup $\bfP\le\bfG$. Let us choose a positive system of roots $\Phi^{+}\subset\Phi(\bfG,\bfT)$ such that $\Phi^{+}\subset\Psi$, and let $\Delta^{+}\subset\Phi^{+}$ be the set of simple positive roots. Let $J\subset\Delta^{+}$ be the subset corresponding to $\Psi$, i.e., $\Psi$ is generated by $\Phi^{+}$ and $-J$.
\begin{defn}[Horizontal degenerations]
Let $\psi\ne\psi_{0}\in(\Lie N)^{*}\bs\lbrace0\rbrace$\footnote{Note that we exclude the trivial character $\rho_{0}$ of $N$. This is because the trivial character of $N$ should be regarded  ``smaller'' than a nontrivial character of $N$, even though both coadjoint orbits are singletons; if we manually define $\dim\cO_{0}:=-\infty$, there can be no horizontal degeneration to the trivial character.}. We call $\psi_{0}$ a \emph{horizontal degeneration} of $\psi$ if the following conditions are satisfied:
\begin{itemize}
\item $\cO_{\psi}\ne\cO_{\psi_{0}}$;
\item $\psi$ and $\psi_{0}$ have the same sizes as $N$-representations (i.e. $\dim\cO_{\psi}=\dim\cO_{\psi_{0}}$);
\item $M\psi_{0}\subset\ov{M\psi}$, given by a one-parameter degeneration in $M$ (i.e., there is a one-parameter subgroup $\lambda:\bG_{m}\rar M$ such that $\psi_{0}=\lim_{t\rar0}\lambda(t)\cdot \psi$);
\item $\lambda$ commutes with $S_{\rho_{\psi}}$.
\end{itemize}
We will denote symbolically as $\psi\rightsquigarrow \psi_{0}$. We call a horizontal degeneration \emph{simple} if the following additional conditions are satisfied:
\begin{itemize}
\item $\dim P\psi=\dim P\psi_{0}+1$;
\item after identifying $\Lie G$ and $(\Lie G)^{*}$ in a $G$-equivariant way via the Killing form, $\delta:=\psi-\psi_{0}$ corresponds to a scalar multiple of a simple negative root orthogonal to $-J$.
\end{itemize}
\end{defn}
\begin{exam}The  degeneration from a non-degenerate character to a nontrivial degenerate character is a simple horizontal degeneration. The degeneration from $\psi_{a,b}$ ($b\ne0$) to $\psi_{a,0}$ as in Example~\ref{exam:GL4} is a simple horizontal degeneration.
\end{exam}
Horizontal degenerations are different from  degenerations arising as the closure relations between coadjoint orbits of $N$, as we require the dimensions of the coadjoint orbits to stay the same along the horizontal degenerations. 

If $\psi\rightsquigarrow\psi_{0}$ is a (simple) horizontal degeneration, then $\rho_{\psi}$ and $\rho_{\psi_{0}}$ seem to be related in many aspects. We record some results in this regard.

\begin{lem}\label{lem:horizontal degeneration S}
Let $\psi\rightsquigarrow\psi_{0}$ be a horizontal degeneration. Then, $S_{\rho_{\psi}}\le S_{\rho_{\psi_{0}}}$.
\end{lem}
\begin{proof}
This follows from the condition that the one-parameter group $\lambda$ commutes with $S_{\rho_{\psi}}$.
\end{proof}

\begin{thm}\label{thm:structure1}
Suppose that $\psi\rightsquigarrow\psi_{0}$ is a simple horizontal degeneration. Suppose that $\wh{S_{\rho_{\psi_{0}}}}$ exists, and choose an extension of $\rho_{\psi_{0}}$ to $\wh{R_{\rho_{\psi_{0}}}}$. Then, $\wh{S_{\rho_{\psi}}}$ exists, and can be given by the pullback of $\wh{S_{\rho_{\psi_{0}}}}$ along the inclusion $S_{\rho_{\psi}}\hrar S_{\rho_{\psi_{0}}}$.
\end{thm}\begin{proof}
By the definition of horizontal degeneration, anything in $S_{\rho_{\psi}}$ clearly stabilizes $\lambda_{0}$, which implies that $S_{\rho_{\psi}}\le S_{\rho_{\psi_{0}}}$. Therefore, $S_{\rho_{\psi}}\times\bG_{m}\le S_{\rho_{\psi_{0}}}$. By (1), we know that $\dim S_{\rho_{\psi_{0}}}=\dim S_{\rho_{\psi}}\times\bG_{m}$. As both are connected (by definition), $S_{\rho_{\psi}}\times\bG_{m}\cong S_{\rho_{\psi_{0}}}$.

Let $Y\in\Lie G$ be a scalar multiple of a simple negative root corresponding to $\delta$ after identifying $\Lie G$ and $(\Lie G)^{*}$ in a  $G$-equivariant way. Let $P_{\psi_{0}}$ be generated by $P$ and $\exp(Y)$. We claim that $\Lie P_{\psi_{0}}=\Lie P\oplus\bC Y$. It suffices to show that $\Lie P\oplus\bC Y$ is a Lie subalgebra of $\Lie G$, or $[Z,Y]\in\Lie P\oplus\bC Y$ for any $Z\in\Lie P$, which is indeed true by (2). 

As $S_{\rho_{\psi}}$ stabilizes both $\psi$ and $\psi_{0}$, it stabilizes $P'$.  Note that $P\bs P_{\psi_{0}}/P=\lbrace 1,w\rbrace$, where $w$ corresponds to the reflection with respect to the simple root parallel to $Y$. As $PwP\subset P_{\psi_{0}}$ is open, by taking the contragredient of \cite[Proposition 1.8]{Bernstein-Zelevinsky-GLn}, we have a $P$-equivariant surjection $\Ind_{N}^{P_{\psi_{0}}}\rho_{\psi}\thrar\Ind^{P}_{wNw^{-1}\cap P}((\rho_{\psi})^{w})$. Note that $wNw^{-1}\cap P=wNw^{-1}\cap N\subset N$. Moreover, if $\alpha\in\Delta^{+}\bs J$ is the positive simple root whose negative corresponds to $Y$, then $s_{\alpha}(\Phi^{+})=\left(\Phi^{+}\bs\lbrace\alpha\rbrace\right)\cup\lbrace-\alpha\rbrace$, so $wNw^{-1}\cap N$ is of codimension one in $N$, and $N=U_{\alpha}\cdot(wNw^{-1}\cap N)$, where $U_{\alpha}$ is the root group corresponding to $\alpha$. 

Let $\fh$ be a polarization subordinate to $\psi$, and let $H=\exp(\fh)$. We may choose so that $\alpha\in \fh$. We claim that $\fh$ is a polarization subordinate to $\psi_{0}$. As $\dim\cO_{\psi}=\dim\cO_{\psi_{0}}$, to show that $\fh$ is a polarization subordinate to $\psi_{0}$, it suffices to show that $\fh$ is a Lie subalgebra subordinate to $\psi_{0}$. Note that, as $wHw^{-1}\cap N\cap U_{\alpha}=1$, $\chi_{\psi}^{w}\rvert_{wHw^{-1}\cap N}=\chi_{\psi_{0}}\rvert_{wHw^{-1}\cap N}$. Therefore, $\fh':=\Lie(wHw^{-1}\cap N)$ is a Lie subalgebra of $\Lie N$ subordinate to $\psi_{0}$. As $\fh'$ is of codimension $1$ in $\fh$, it suffices to check that $\psi_{0}([\alpha,X])=0$ for any $X\in\fh'$. As $\fh$ is already a polarization subordinate to $\psi$, $\psi([\alpha,X])=0$. Thus, it suffices to show that $\delta([\alpha,X])=0$. As the root space decomposition of $[\alpha,X]$ only involves roots bigger than $\alpha$, $\psi_{0}([\alpha,X])=0$, as desired.

Therefore, 
\[(\rho_{\psi})^{w}\rvert_{wNw^{-1}\cap N}=\left(\ind_{wHw^{-1}}^{wNw^{-1}}\chi_{\psi}^{w}\right)\rvert_{wNw^{-1}\cap N}=\ind^{wNw^{-1}\cap N}_{wHw^{-1}\cap N}\chi_{\psi}^{w}=\ind^{wNw^{-1}\cap N}_{wHw^{-1}\cap N}\chi_{\psi_{0}},\]
\[\rho_{\psi_{0}}\rvert_{wNw^{-1}\cap N}=\left(\ind_{H}^{N}\chi_{\psi_{0}}\right)\rvert_{wNw^{-1}\cap N}=\ind_{wNw^{-1}\cap H}^{wNw^{-1}\cap N}\chi_{\psi_{0}}.\]
Note that $wHw^{-1}\cap H= wHw^{-1}\cap N=wNw^{-1}\cap H$. Therefore, both $ (\rho_{\psi})^{w}\rvert_{wNw^{-1}\cap N}$ and $\rho_{\psi_{0}}\rvert_{wNw^{-1}\cap N}$ are realized in the same representation space. Note that $S_{\rho_{\psi}}$ stabilizes $\alpha$, so $\wh{S_{\rho_{\psi_{0}}}}\ltimes N$ has a subgroup $\wh{S_{\rho_{\psi}}}\ltimes(wNw^{-1}\cap N)$, where $\wh{S_{\rho_{\psi}}}$ is defined by the pullback central cover along the natural inclusion $S_{\rho_{\psi}}\hrar S_{\rho_{\psi}}\times\bG_{m}\cong S_{\rho_{\psi_{0}}}$. Therefore, the chosen extension of $\rho_{\psi_{0}}$ to $\wh{S_{\rho_{\psi_{0}}}}\ltimes N$ gives rise to the extension of $(\rho_{\psi})^{w}\rvert_{wNw^{-1}\cap N}$ to $\wh{S_{\rho_{\psi}}}\ltimes(wNw^{-1}\cap N)$, or the extension of $\rho_{\psi}\rvert_{N\cap w^{-1}Nw}$ to $w^{-1}\wh{S_{\rho_{\psi}}}w\ltimes (N\cap w^{-1}Nw)$. As $S_{\rho_{\psi}}$ stabilizes $\alpha$, $w^{-1}\wh{S_{\rho_{\psi}}}w\cong \wh{S_{\rho_{\psi}}}$. There are two actions on $\rho_{\psi}$, the action by $N$ and the action by $\wh{S_{\rho_{\psi}}}\ltimes (N\cap w^{-1}Nw)$. These two actions give rise to an action of $\wh{S_{\rho_{\psi}}}\ltimes N$ as $\wh{S_{\rho_{\psi}}}$ stabilizes $U_{\alpha}$.
\end{proof}
The proof of Theorem~\ref{thm:structure1} implies that there is a natural $P$-equivariant surjection
\[\Ind_{N}^{P_{\psi_{0}}}\rho_{\psi}\thrar\Ind_{wNw^{-1}\cap N}^{P}\rho_{\psi_{0}}.\]
The following example suggests a relation between $\DWh_{\psi}(\pi)$ and $\DWh_{\psi_{0}}(\pi)$.
\begin{exam}

As noted in Remark~\ref{rem:GL4}, if we use the notations of Example~\ref{exam:GL4}, for $a,b\in F^{\times}$,  the generalized Whittaker models associated with $\psi_{a,b}$ is infinite-dimensional, whereas Theorem~\ref{thm:GLn} implies that the generalized Whittaker models associated with $\psi_{a,0}$ is in general of multiplicity one. 

Let $\pi$ be an irreducible supercuspidal reprsentation of $\GL_{4}(F)$. Then,
\[\DWh_{\psi_{a,b}}(\pi)=\Hom_{H}(\pi,\chi_{\psi_{a,b}}),\quad \DWh_{\psi_{a,0}}(\pi)=\Hom_{H}(\pi,\chi_{\psi_{a,0}}).\]
Let\[H_{1}=\left\lbrace\begin{pmatrix}1& 0 & 0 & *\\
0 & 1 & 0 & *\\
0 & 0 & 1 & *\\
0 & 0 & 0 & 1\end{pmatrix}\right\rbrace,\quad H_{2}=\left\lbrace\begin{pmatrix} 1& 0 & 0 & 0\\
0 & 1 & * & 0\\
0 & 0 & 1 & 0\\
0 & 0 & 0 & 1\end{pmatrix}\right\rbrace,\quad B_{2}=\left\lbrace\begin{pmatrix} 1& 0 & 0 & 0\\
0 & * & * & 0\\
0 & 0 & * & 0\\
0 & 0 & 0 & 1\end{pmatrix}\right\rbrace,\]\[ U_{3}=\left\lbrace\begin{pmatrix}1 & 0 & 0 & 0\\
* & 1 & 0 & 0\\
* & 0 & 1 & 0\\
0 & 0 & 0 & 1\end{pmatrix}\right\rbrace,\quad P_{2}=\left\lbrace\begin{pmatrix} 1 & 0 & 0 & 0\\
0& 1 & 0 & 0\\ 0 & * & * & 0\\ 0 & 0 & 0 & 1\end{pmatrix}\right\rbrace,\quad P_{3}=\left\lbrace\begin{pmatrix} 1 & 0 & 0 & 0\\
* & * & * & 0\\ * & * & * & 0\\ 0 & 0 & 0 & 1\end{pmatrix}\right\rbrace.\]
Let $\chi_{1}=\chi_{a,0}\rvert_{H_{1}}=\chi_{a,b}\rvert_{H_1}$ and $\chi_{2}=\chi_{a,b}\rvert_{H_{2}}$. Then, as $H_{1}$ is a normal subgroup of $H$, 
\[\DWh_{\psi_{a,b}}(\pi)=\Hom_{H_{2}}(J_{H_{1},\chi_{1}}(\pi),\chi_{2}),\quad \DWh_{\psi_{a,0}}(\pi)=\Hom_{H_{2}}(J_{H_{1},\chi_{1}}(\pi),\bfone).\]By the same reason as the proof of Theorem~\ref{thm:GLn}, $J_{H_{1},\chi_{1}}(\pi)=\ind_{P_{2}U_{3}}^{P_{3}}(\tau\otimes\psi)$, where $\tau$ is an irreducible smooth admissible non-degenerate representation of $P_{2}$ trivially exteded over $U_{3}$, and $\psi:U_{3}\rar\bC^{\times}$ is defined as
\[\psi\begin{pmatrix}1&0&0&0
\\
x_{2,1}&1&0&0\\
x_{3,1}&0&1&0\\
0&0&0&1\end{pmatrix}:=\varepsilon(x_{2,1}).\]
By  Mackey theory, $J_{H_{1},\chi_{1}}(\pi)$ sits in a short exact sequence
\[0\rar\ind_{P_{2}U_{3}\cap B_{2}}^{B_{2}}(\tau\otimes\psi)\rar J_{H_{1},\chi_{1}}(\pi)\rar\ind_{wP_{2}U_{3}w^{-1}\cap B_{2}}^{B_{2}}((\tau\otimes\psi)^{w})\rar0,\]
where $w=\begin{pmatrix}1&0&0&0\\0&0&-1&0\\0&1&0&0\\0&0&0&1\end{pmatrix}$. 
Since $\psi$ is trivial on $B_{2}$, we have
\[0\rar\ind_{P_{2}U_{3}\cap B_{2}}^{B_{2}}(\tau)\rar J_{H_{1},\chi_{1}}(\pi)\rar\ind_{wP_{2}U_{3}w^{-1}\cap B_{2}}^{B_{2}}(\tau^{w})\rar0.\]
Since $\chi_{2}$ and $\bfone$ are  injective objects in $\Rep^{\sm}(H_{2})$ (which is a consequence of the exactness of the usual Jacquet functor), we have
\[0\rar\Hom_{H_{2}}(\ind_{wP_{2}U_{3}w^{-1}\cap B_{2}}^{B_{2}}(\tau^{w}),\chi_{2})\rar\DWh_{\psi_{a,b}}(\pi)\rar\Hom_{H_{2}}(\ind_{P_{2}U_{3}\cap B_{2}}^{B_{2}}(\tau),\chi_{2})\rar0.\]
\[0\rar\Hom_{H_{2}}(\ind_{wP_{2}U_{3}w^{-1}\cap B_{2}}^{B_{2}}(\tau^{w}),\bfone{})\rar\DWh_{\psi_{a,0}}(\pi)\rar\Hom_{H_{2}}(\ind_{P_{2}U_{3}\cap B_{2}}^{B_{2}}(\tau),\bfone{})\rar0.\]
Note that the exact sequences are $S_{a,b}':=\lbrace\diag(1,x,x,1)~:~x\in F^{\times}\rbrace$-equivariant and $S_{a,0}':=\lbrace\diag(1,x,y,1)~:~x,y\in F^{\times}\rbrace$-equivariant, respectively. 

Note that the proof of Theorem~\ref{thm:GLn} shows that $\Hom_{H_{2}}(\ind_{wP_{2}U_{3}w^{-1}\cap B_{2}}^{B_{2}}(\tau^{w}),\bfone)=0$, and for any smooth irreducible admissible representation $\sigma$ of $S_{a,0}'$, $\Hom_{S_{a,0}'}(\sigma^{\vee},\Hom_{H_{2}}(\ind_{P_{2}U_{3}\cap B_{2}}^{B_{2}}(\tau),\bfone))$ is one-dimensional. On the other hand, it is easy to see that, for any smooth irreducible admissible representation $\sigma$ of $S'_{a,b}$, 
\[\dim\Hom_{S'_{a,b}}(\sigma^{\vee},\Hom_{H_{2}}(\ind_{wP_{2}U_{3}w^{-1}\cap B_{2}}^{B_{2}}(\tau^{w}),\chi_{2}))=1,\] as $\tau^{w}$ is a smooth irreducible admissible non-degenerate representation of the $2\times 2$ mirabolic group $wP_{2}U_{3}w^{-1}\cap B_{2}$, whereas\[\dim\Hom_{S'_{a,b}}(\sigma^{\vee},\Hom_{H_{2}}(\ind_{P_{2}U_{3}\cap B_{2}}^{B_{2}}(\tau),\chi_{2}))=\infty, \]as this space contains $\Hom_{S'_{a,b}H_{2}}(\chi_{2}^{-1}\otimes\sigma,\Ind^{B_{2}}_{\lbrace1\rbrace}(\bfone))$, which is infinite-dimensional, as a subspace; this follows from the fact that $\tau\rvert_{P_{2}U_{3}\cap B_{@}}\cong\ind_{\lbrace1\rbrace}^{P_{2}U_{3}\cap B_{2}}(\bfone)$, and the $H_{2}$-contragredient of $\ind_{P_{2}U_{3}\cap B_{2}}^{B_{2}}(\ind_{\lbrace1\rbrace}^{P_{2}U_{3}\cap B_{2}}(\bfone))$ contains the $B_{2}$-contragredient of $\ind_{P_{2}U_{3}\cap B_{2}}^{B_{2}}(\ind_{\lbrace1\rbrace}^{P_{2}U_{3}\cap B_{2}}(\bfone))$, which is $\Ind_{P_{2}U_{3}\cap B_{2}}^{B_{2}}(\Ind_{\lbrace1\rbrace}^{P_{2}U_{3}\cap B_{2}}(\bfone))$. 

The above calculations suggest that a part of the degenerate Whittaker model has the corresponding generalized Whittaker model multiplicity-free: for $\DWh_{\psi_{a,b}}(\pi)$, it is the natural subspace $\Hom_{H_{2}}(\ind_{wP_{2}U_{3}w^{-1}\cap B_{2}}^{B_{2}}(\tau^{w}),\chi_{2})$; for $\DWh_{\psi_{a,0}}(\pi)$, it is the natural quotient $\Hom_{H_{2}}(\ind_{P_{2}U_{3}\cap B_{2}}^{B_{2}}(\tau),\bfone)$. The change of multiplicities of the two parts of $\DWh$ under the simple horizontal degeneration $\psi_{a,b}\rightsquigarrow\psi_{a,0}$ is summarized below as follows.
\begin{center}
\begin{tabular}{ | c | c | c | c | }
\hline
  & $\psi_{a,b}$ & $\rightsquigarrow$ & $\psi_{a,0}$ \\ 
\hline
 The ``sub'' of $\DWh$ (corresponding to $\ind_{wP_{2}U_{3}w^{-1}\cap B_{2}}^{B_{2}}(\tau^{w})$) & $1$ & $\rightsquigarrow$ & $0$ \\  
\hline
 The ``quo'' of $\DWh$ (corresponding to $\ind_{P_{2}U_{3}\cap B_{2}}^{B_{2}}(\tau)$) & $\infty$ & $\rightsquigarrow$ & $1$\\
\hline
\end{tabular}
\end{center}

It will be desirable to determine, for more general horizontal degenerations $\psi\rightsquigarrow\psi_{0}$, the subspace of finite multiplicity of $\DWh_{\psi}(\pi)$ based on $\DWh_{\psi_{0}}(\pi)$.
\end{exam}
\section{{Further questions}}
We raise further questions regarding the generalized and degenerate Whittaker models.
\subsection{Hamiltonian spaces and quantization}\label{sec:BZSV}
Kirillov's orbit method suggests that representations should be constructed as the quantization of coadjoint orbits. Put further, a recent program of \cite{BZSV} suggests to study the branching problems in terms of quantization of Hamiltonian varieties. In \cite{GanPeng}, the generalized Whittaker models in the classical sense (the generalized B-Whittaker models as in Remark~\ref{rem:history}) were studied in this optic. We briefly explain how the generalized Whittaker models of Definition~\ref{def:GWh} arise as the quantization of Hamiltonian varieties, and compare them with the above literatures. We will freely use the terminologies used in them.

For $\psi\in(\Lie N)^{*}$, the coadjoint orbit $\cO_{\psi}$ is a Hamiltonian $N$-variety; it is well-known that $\cO_{\psi}$ is equipped with the Kirillov--Kostant--Souriau symplectic form, which is $N$-invariant, and the moment map is given by the inclusion $\cO_{\psi}\rar(\Lie N)^{*}$. Note that the inductive procedure of Kirillov in proving Theorem~\ref{thm:Kirillov} can also be regarded as realizing the coadjoint orbit as a \emph{twisted cotangent bundle}. More precisely, if we take a polarization $\fh\subset\Lie N$ subordinate to $\psi$, then \[\cO_{\psi}\cong(\psi+\fh^{\perp})\times^{H}N,\]where $\psi+\fh^{\perp}$ is the preimage of the image of $\psi$ along the natural surjection $(\Lie N)^{*}\thrar\fh^{*}$, and $H=\exp(\fh)$. This is the symplectic induction of the Hamiltonian $H$-variety $\psi$, which is a singleton with a moment map sending the point to $\psi\rvert_{\fh}$. 

As $S_{\rho_{\psi}}$ preserves the Kirillov--Kostant--Souriau form of $\cO_{\psi}$, $\cO_{\psi}$ can be regarded as a Hamiltonian $R_{\rho_{\psi}}$-space. For  a Hamiltonian $S_{\rho_{\psi}}$-space $Y$, we can also regard it as a Hamiltonian $R_{\rho_{\psi}}$-space. Then, the symplectic induction of $Y\times \cO_{\psi}$ from $R_{\rho_{\psi}}$ to $G$,
\[X_{Y,\psi}:=(Y\times \cO_{\psi})\times^{R_{\rho_{\psi}}}_{(\Lie R_{\rho_{\psi}})^{*}}T^{*}G\cong((Y\times \cO_{\psi})\times_{(\Lie R_{\rho_{\psi}})^{*}}(\Lie G)^{*})\times^{R_{\rho_{\psi}}}G,\]
is a Hamiltonian $G$-space that should govern the generalized Whittaker model $\GWh_{\sigma,\psi}$ for $\sigma$ ``corresponding to'' $Y$ (the ``anomaly'' should also be considered when quantizing it). As the moment map $\cO_{\psi}\rar(\Lie R_{\rho_{\psi}})^{*}$ is injective, we can rewrite the above expression as in \cite[4.2.1]{GanPeng}: namely,\[(Y\times \cO_{\psi})\times_{(\Lie R_{\rho_{\psi}})^{*}}(\Lie G)^{*}\cong Y\times_{(\Lie S_{\rho_{\psi}})^{*}}(\cO_{\psi}+(\Lie N)^{\perp}),\]
where $\cO_{\psi}+(\Lie N)^{\perp}\subset(\Lie G)^{*}$ is the preimage of $\cO_{\psi}$ along the natural surjection $(\Lie G)^{*}\thrar(\Lie N)^{*}$.  

Recall that, in \cite{BZSV}, a Hamiltonian $G$-variety $M$ is called \emph{hyperspherical} if, in addition to the multiplicity-free condition, several auxiliary conditions are satisfied; in particular, $M$ is always required to be smooth, affine, and ``neutral.'' In general, these auxiliary  conditions fail to hold for the Hamiltonian $G$-space $X_{Y,\psi}$. For example, without the $\fsl_{2}$-triple structure, the neutrality will not hold in general.
\begin{exam}
Consider the setup of Example~\ref{exam:GL4}. By Theorem~\ref{thm:GLn}, one expects that $X_{0,\psi_{a,0}}$ is close to being multiplicty-free, where $a\in F^{\times}$ and $0$ denotes the trivial Hamiltonian $S_{\rho_{a,0}}$-space. Note that
\[X_{0,\psi_{a,0}}=((\cO_{\psi_{a,0}}+(\Lie N)^{\perp})\cap(\Lie S_{\rho_{a,0}})^{\perp})\times^{R_{\rho_{\psi_{a,0}}}}G,\]where $(\Lie S_{\rho_{a,0}})^{\perp}\subset(\Lie G)^{*}$ is the kernel of the natural surjection $(\Lie G)^{*}\thrar(\Lie S_{\rho_{a,0}})^{*}$. One can check directly that $(\cO_{\psi_{a,0}}+(\Lie N)^{\perp})\cap(\Lie S_{\rho_{a,0}})^{\perp}$ is a \emph{quadric hypersurface} of $(\Lie S_{\rho_{a,0}})^{\perp}$, so $X_{0,\psi_{a,0}}\rar G/R_{\rho_{\psi_{a,0}}}$ is affine, where $G/R_{\rho_{\psi_{a,0}}}$ is $\bG_{m}^{2}$-bundle over the flag variety of $\GL_{4}$ over $F$. Therefore, $X_{0,\psi_{a,0}}$ is \emph{not affine}.
\end{exam}
On the other hand, as observed before, there are certainly many generalized Whittaker models of multiplicity one that are new. One may expect that these generalized Whittaker models (or rather their corresponding Hamiltonian spaces) fit into a more generalized version of the relative Langlands duality of \cite{BZSV} where the auxiliary conditions as mentioned above are relaxed. 
\subsection{Complete Fourier expansion of automorphic forms}\label{sec:Fourier}
For this subsection only, let $\bfG$ be a semisimple  algebraic group over $\bQ$. 
Let $\pi=\bigotimes_{p\le \infty}\pi_{p}$ is a cuspidal automorphic representation of $\bfG(\bA)$, where $\pi_{p}$ is a smooth admissible irreducible $\bfG(\bQ_{p})$-representation for a finite prime $p<\infty$, and $\pi_{\infty}$ is a $(\fg,K_{\infty})$-module, where $\fg=\Lie \bfG(\bC)$ and $K_{\infty}\le \bfG(\bR)$ is a maximal compact subgroup. 
Let  $f:\bfG(\bA)\rar\bC$ be a cuspidal automorphic form in $\pi$.

Let $\bfP\le \bfG$ be a parabolic group and $\bfP=\bfM\bfN$ be the Levi decomposition. In many classical contexts, such as the $q$-expansion of a (Hilbert) modular form or the Fourier--Jacobi expansion of a Siegel modular form, it is very useful to take the Fourier expansion of $f$ along $\bfN$. More concretely, for $g\in \bfG(\bA)$, we are interested in describing the map $f_{g}:\bfN(\bA)\rar\bC$, $n\mapsto f(ng)$. Note that this map is left-$\bfN(\bQ)$-invariant, so defines $f_{g}\in L^{2}(\bfN(\bQ)\bs \bfN(\bA))$. 
Note that the Kirillov's orbit method works for $\bfN(\bQ_{p})$ for all primes $p\le \infty$ (including the archimedean place), and all irreducible representations of $\bfN(\bQ_{p})$ are unitary. Let $\fn=\Lie \bfN$ be the Lie algebra of $\bfN$ over $\bQ$. Then, the description of the Plancherel measure of $\bfN(\bQ_{p})$ gives the following description,
\[L^{2}(\bfN(\bQ)\bs \bfN(\bA))=\wh{\bigoplus}_{\psi=(\psi_{p})\in\left(\prod'_{p\le\infty}(\fn(\bQ_{p})^{*}/\bfN(\bQ_{p}))\right)^{\bfN(\bQ)}}\left(\bigotimes'_{p\le\infty}\rho_{\psi_{p}}\right)\otimes W_{\psi},\]where $W_{\psi}$ is a multiplicity space. Let us elaborate more on the index set of the Hilbert space direct sum. For all but finitely many  $p$, there exists $\bfN(\bZ_{p})$, and an irreducible smooth admissible representation $\rho$ of $\bfN(\bQ_{p})$ is spherical if $\rho$ has a nonzero $\bfN(\bZ_{p})$-invariant vector. It is easy to see that, by Mackey theory (e.g., \cite[5.5]{Vigneras-book}), for $\psi_{p}\in\fn(\bQ_{p})^{*}$, $\rho_{\psi_{p}}$ is spherical if and only if there exists a polarization $\fh\subset\fn(\bQ_{p})$ subordinate to $\psi_{p}$ such that $\ker\psi_{p}$ contains $\log(\bfN(\bZ_{p}))\cap \fh$. We define $\prod'_{p\le\infty}(\fn(\bQ_{p})^{*}/\bfN(\bQ_{p}))$ to be the subset of $\prod_{p\le\infty}(\fn(\bQ_{p})^{*}/\bfN(\bQ_{p}))$ where all but finitely many entries are spherical. This set admits a conjugation action by $\bfN(\bQ)$, and the Hilbert space direct sum runs over the $\bfN(\bQ)$-fixed elements. According to this decomposition, we have an infinite sum
\[f_{g}(n)=\sum_{\psi\in\left(\prod'_{p\le\infty}(\fn(\bQ_{p})^{*}/\bfN(\bQ_{p}))\right)^{\bfN(\bQ)}}f_{g,\psi}(n).\]
For any linear functional $\lambda_{\psi}:W_{\psi}\rar\bC$, the function $g\in\bfG(\bA)\mapsto\lambda_{\psi}(f_{g,\psi})$ defines a vector in $\bigotimes'_{p\le \infty}\Ind_{\bfN(\bQ_{p})}^{\bfG(\bQ_{p})}(\rho_{\psi_{p}})$. Varying $f$, we obtain a vector in $\prod_{p\le\infty}\DWh_{\psi_{p}}(\pi_{p})$. One may  refine this further by using the action of the metaplectic cover of the stabilizer group $\wh{S_{\rho_{\psi_{p}}}}$ on $W_{\psi}$; decomposing further using this action will give us a vector in the product of generalized Whittaker models. This does not in general give an infinite sum expression of an automorphic form as one may want, due to the existence of continuous spectrum of the stabilizer groups;  if $\bfM(\bQ_{p})$ is compact modulo center for all $p\le\infty$ (e.g., when $M$ is abelian), however, there is no such a problem. 

A Fourier expansion of an automorphic form, given the explicit realizations of the generalized / degenerate Whittaker model spaces, will attach the Fourier \emph{coefficients} to an automorphic form as numbers. As in the case of $q$-expansion of modular forms, it will be interesting to study the algebraicity properties of such numbers. For example, Pollack has shown the arithmeticity of Fourier coefficients of quaternionic modular forms, e.g., \cite{Pollack1,Pollack2,Pollack3}. See also \cite{Narita1,Narita2,Narita3}.

\appendix
\section{{Smooth representation theory}}\label{sec:appendixA}
We present a series of lemmas in smooth representation theory of td-groups, generally regarding the functor of taking smooth parts $V\mapsto V_{\infty}$.  Throughout this section, $G$ will be a td-group (see {Notations and Conventions}).

\begin{lem}[Universal property of smooth part]\label{lem:univprop}
Let $(\pi,V)$ be a representation of $G$. Then, the natural inclusion $\pi_{\infty}\rar \pi$ of the smooth part of $\pi$ is universal among any $G$-linear maps $\sigma\rar \pi$ from a smooth representation $\sigma$ of $G$ to $\pi$.
\end{lem}
\begin{proof}
Since any $G$-linear map sends a smooth vector to a smooth vector, given any $G$-linear map $\sigma\rar \pi$ from a smooth representation $\sigma$, it must factor through $\sigma\rar\pi_{\infty}\rar\pi$. This factorization is unique as $\pi_{\infty}\rar\pi$ is injective.
\end{proof}

\begin{coro}\label{coro:smoothHom}
Given a representation $(\pi,V)$ of $G$ and a smooth representation $(\sigma,W)$, the natural inclusion $\pi_{\infty}\rar\pi$ induces an isomorphism
\[\Hom_{G}(\sigma,\pi_{\infty})\riso\Hom_{G}(\sigma,\pi).\]
\end{coro}
\begin{proof}
Obvious from Lemma~\ref{lem:univprop}.
\end{proof}

\begin{coro}[{\cite[p. 27]{Prasad-ICM}}]\label{coro:V1V2V3}
If $(\pi_{1},V_{1})$, $(\pi_{2},V_{2})$, $(\pi_{3},V_{3})$ are smooth representations of $G$, then
\[\Hom_{G}(V_{1}\otimes V_{2},V_{3}^{\vee})\cong\Hom_{G}(V_{1}\otimes V_{3},V_{2}^{\vee}).\]
\end{coro}
\begin{proof}
This follows from the usual adjunction of $G$-linear maps
\[\Hom_{G}(V_{1}\otimes V_{2},V_{3}^{*})\cong \Hom_{G}(V_{1}\otimes V_{3},V_{2}^{*}),\]and Corollary~\ref{coro:smoothHom}.
\end{proof}

\begin{defn}
Let $K$ be an open compact subgroup of $G$, and let $(\pi,V)$ be a smooth representation of $G$. It is well-known that $V$ is $K$-semisimple (e.g. \cite[Lemma 1.2.2]{BushnellHenniart}). For $\rho\in\wh{K}$, let $V(\rho)$ be the $\rho$-isotypic part of $V$.
\end{defn}

\begin{prop}\label{prop:K-isotypic}
Let $K$ be an open compact subgroup of $G$, and let $(\pi,V)$ be a smooth admissible representation of $G$. Let $\rho\in\wh{K}$, and let $\rho^{*}$ be its dual representation. Then, we have an isomorphism
\[\left(V^{\vee}\right)(\rho)\cong\left(V(\rho^{*})\right)^{*}.\]
\end{prop}

\begin{proof}
Obvious from the description $V^{\vee}=\bigoplus_{\rho\in\wh{K}}\left(V(\rho)\right)^{*}$ (e.g. \cite[Chapter 4, (2.7)]{Bump}).
\end{proof}

\begin{coro}\label{coro:contragredient-tensor}
Let $(\pi_{1},V_{1})$ and $(\pi_{2},V_{2})$ be smooth admissible $G$-representations such that $(\pi_{1}\otimes \pi_{2},V_{1}\otimes V_{2})$ is also admissible as a $G$-representation\footnote{This assumption is necessary. Namely, it is \emph{not true} in general that $(V_{1}\otimes V_{2})^{\vee}\cong V_{1}^{\vee}\otimes V_{2}^{\vee}$ for smooth admissible $G$-representations $V_{1},V_{2}$; otherwise, $\left((V_{1}\otimes V_{2})^{\vee}\right)^{\vee}\cong\left(V_{1}^{\vee}\otimes V_{2}^{\vee}\right)^{\vee}\cong \left(V_{1}^{\vee}\right)^{\vee}\otimes \left(V_{2}^{\vee}\right)^{\vee}\cong V_{1}\otimes V_{2}$, which implies that $V_{1}\otimes V_{2}$ is admissible, and this is not in general true.}. Then,
\[(V_{1}\otimes V_{2})^{\vee}\cong V_{1}^{\vee}\otimes V_{2}^{\vee}.\]
\end{coro}

\begin{proof}
Note that, under the natural injection $V_{1}^{*}\otimes V_{2}^{*}\hrar(V_{1}\otimes V_{2})^{*}$, a tensor product of smooth vectors is smooth, so there is a natural injection $V_{1}^{\vee}\otimes V_{2}^{\vee}\hrar (V_{1}\otimes V_{2})^{\vee}$. We would like to show that this is an isomorphism.

Let $K\le G$ be an open compact subgroup. Then, it suffices to show that \[\tag{\textreferencemark
}\label{eq:A3}
\left(V_{1}^{\vee}\otimes V_{2}^{\vee}\right)^{K}\hrar \left((V_{1}\otimes V_{2})^{\vee}\right)^{K}\cong\left((V_{1}\otimes V_{2})^{K}\right)^{*}\] is surjective. By Proposition~\ref{prop:K-isotypic}, (\ref{eq:A3}) is
\[\left(\left(\bigoplus_{\rho\in\wh{K}}\left(V_{1}(\rho)\right)^{*}\right)\otimes\left(\bigoplus_{\tau\in\wh{K}}\left(V_{2}(\tau)\right)^{*}\right)\right)^{K}\hrar\left(\left(\left(\bigoplus_{\rho\in\wh{K}}V_{1}(\rho)\right)\otimes\left(\bigoplus_{\tau\in\wh{K}}V_{2}(\tau)\right)\right)^{K}\right)^{*}.\]By the assumption that $V_{1}\otimes V_{2}$ is admissible, there are only finitely many pairs $(\rho,\tau)$ of smooth irreducible representations of $K$ such that $V_{1}(\rho)\ne0$, $V_{2}(\tau)\ne0$ and $(\rho\otimes\tau)^{K}\ne0$. Furthermore, for such a pair $(\rho,\tau)$, $\dim V_{1}(\rho)<\infty$ and $\dim V_{2}(\tau)<\infty$ as $\rho$ and $\tau$ are finite-dimensional, and $\dim(\rho\otimes\tau)^{K}=\dim(\rho^{*}\otimes\tau^{*})^{K}$ as every finite-dimensional smooth representation of $K$ is semisimple. Therefore, if we enumerate these pairs as $(\rho_{i},\tau_{i})$, $i=1,\cdots,n$, then (\ref{eq:A3}) is\[\left(\bigoplus_{i=1}^{n}\left(V_{1}(\rho_{i})\right)^{*}\otimes \left(V_{2}(\tau_{i})\right)^{*})\right)^{K}\hrar\left(\left(\bigoplus_{i=1}^{n}V_{1}(\rho_{i})\otimes V_{2}(\tau_{i})\right)^{K}\right)^{*},\]
which is an injection between finite-dimensional vector spaces of the same dimension, so it is an isomorphism.
\end{proof}

\begin{coro}\label{coro:random}
Let $G=H\ltimes N$ be a semidirect product of td-groups $H,N$. Let $(\pi,V)$ be a smooth $G$-repersentation. Let $(\rho,X)$ be a smooth admissible $G$-representation, and $(\sigma,Y)$ be a smooth admissible $H$-representation, seen as a smooth admissible $G$-representation via pullback. Suppose that $(\rho\otimes \sigma,X\otimes Y)$ is admissible as a $G$-representation. Then, 
\[\Hom_{G}(V,X\otimes Y)\cong \Hom_{H}(Y^{\vee},\Hom_{N}(V,X))\cong\Hom_{H}(Y^{\vee},\Hom_{N}(V,X)_{\infty}),\]where $\Hom_{N}(V,X)_{\infty}$ is the $H$-smooth part of $\Hom_{N}(V,X)$.
\end{coro}
\begin{proof}
We have
\begin{align*}
\Hom_{G}(V,X\otimes Y)
&
\cong\Hom_{G}(V\otimes(X\otimes Y)^{\vee},\bfone)&\text{ (Corollary~\ref{coro:V1V2V3}; $((X\otimes Y)^{\vee})^{\vee}=X\otimes Y$)}
\\ & 
\cong\Hom_{G}(V\otimes X^{\vee}\otimes Y^{\vee},\bfone)&\text{ (Corollary~\ref{coro:contragredient-tensor})}
\\ &
\cong\Hom_{G}(V\otimes Y^{\vee},X)&\text{ (Corollary~\ref{coro:V1V2V3}; $(X^{\vee})^{\vee}=X$)}
\\ &
\cong\Hom_{G}(Y^{\vee},\Hom(V,X))&\text{ (Hom-tensor adjunction)}
\\ &
\cong\Hom_{H}(Y^{\vee},\Hom_{N}(V,X))
\\ &
\cong\Hom_{H}(Y^{\vee},\Hom_{N}(V,X)_{\infty})&\text{ (Corollary~\ref{coro:smoothHom}),}
\end{align*}
where $\Hom(V,X)$ is the internal Hom, a $G$-representation defined by the formula (\ref{eq:Homrep}). Note that there is no ambiguity of the meaning of $Y^{\vee}$ as the contragredient as an $H$-representation is the same as the contragredient as a $G$-representation.
\end{proof}

\bibliographystyle{alpha}
\bibliography{bib}
\end{document}